\documentstyle[12pt,fleqn,dina4,amssymb]{article}

\def\x{\times}
\def\ds{\displaystyle}

\def\eset{\emptyset}

\def\bigspitem{\par\hangone\textputone}
\def\hangone{\hangindent 3\parindent}
\def\textputone#1{\noindent
	 \hbox to 3\parindent{\hss#1\enspace}\ignorespaces}

\begin{document}
\sloppy

\newcommand{\fl}{\mbox{ {\rm fl}} }
\newcommand{\afl}{\mbox{ {\rm afl}} }
\newcommand{\acl}{\mbox{ {\rm acl}} }
\newcommand{\hcl}{\mbox{ {\rm hcl}} }
\newcommand{\gsh}{\mbox{ {\rm gsh}} }
\newcommand{\sh}{\mbox{ {\rm sh}} }
\newcommand{\GG}{\mathop{\rm I\!\!\! G}\nolimits}
\newcommand{\Id}{\mbox{ {\rm Id}} }
\newcommand{\id}{\mbox{ {\rm id}} }
\newcommand{\gf}{\mbox{ {\rm gf}} }
\newcommand{\cl}{\mbox{ {\rm cl}} }
\newcommand{\dpo}{\mathop{ \to\!\!\!\!\!\to }\nolimits}
\newcommand{\supp}{\mbox{ {\rm supp}} }
\newcommand{\ver}{\mbox{ {\rm ver}} }
\newcommand{\edg}{\mbox{ {\rm edg}} }
\newcommand{\diam}{\mbox{ {\rm diam}} }
\newcommand{\FP}{\mbox{ {\rm FP}} }
\newcommand{\FH}{\mbox{ {\rm FH}} }
\newcommand{\Res}{\mbox{ {\rm Res}} }
\newcommand{\dpm}{\mathop{\ \succ\!\!\to\!\!\!\!\!\to }\nolimits}
\newcommand{\pfm}{\mathop{\ \succ\!\!\to}\nolimits}
\newcommand{\NN}{\mathop{\rm I\! N}\nolimits}
\newcommand{\QQ}{\mathop{\rm I\!\!\! Q}\nolimits}
\newcommand{\RR}{\mathop{\rm I\! R}\nolimits}
\newcommand{\ZZ}{\mathop{\sl Z\!\! Z}\nolimits}
\newcommand{\BB}{\mathop{\rm I\! B}\nolimits}
\newcommand{\CC}{\mathop{\rm C\!\!\! I}\nolimits}
\newcommand{\FF}{\mathop{\rm I\! F}\nolimits}
\newcommand{\Hom}{\mathop{\rm Hom}\nolimits}
\newcommand{\HNN}{\mathop{\rm Hnn}\nolimits}
\newcommand{\Isom}{\mathop{\rm Isom}\nolimits}
\newcommand{\EE}{\mathop{\rm I\! E}\nolimits}
\newcommand{\HH}{\mathop{\rm I\! H}\nolimits}
\newcommand{\SL}{\mathop{ {\rm SL}} }
\newcommand{\GL}{\mathop{ {\rm GL}} }
\renewcommand{\thefootnote}{\fnsymbol{footnote}}

\centerline{{\Large \bf Connectivity properties of group}}
\centerline{{\Large \bf actions on non-positively curved spaces II:}}
\centerline{{\Large \bf The geometric invariants}}
~\\

\centerline{\large by}
~\\

\centerline{\large Robert Bieri\footnote{The first-named author was
supported in part by a grant from the Deutsche Forschungsgemeinschaft.}
and Ross Geoghegan\footnote{The second-named
author was supported in part by a grant from the National Science Foundation.}}
~\\

\section*{Introduction}
\renewcommand{\thefootnote}{\arabic{footnote}}
\setcounter{footnote}{20}

This Part II continues our paper [BG$_{\mbox{I}}$], but we have made it as
independent of that paper as possible.  The reader familiar with the
self-contained essay on finitary sheaves and finitary maps which is \S4 of
[BG$_{\mbox{I}}$] can read this second paper with occasional references back
to [BG$_{\mbox{I}}$] for some details and analogies.  A full outline of the
present paper is found in \S10 below.

Let $G$ be a group, let $M$ be a simply connected ``non-positively curved'',
i.e. CAT$(0)$, metric space, and let $\rho : G \to \mbox{Isom}(G)$ be an
action of $G$ on $M$ by isometries.  In our previous paper [BG$_{\mbox{I}}$]
we introduced, when $G$ is of type $F_n$, a property of the action $\rho$
which we called ``controlled $(n-1)$-connectedness'', abbreviated $CC^{n-1}$.
This property was defined in terms of the filtration of $M$ by the balls
$B_r(a)$ centered at some base point $a\in M$ together with a free contractible 
$G$-CW-complex over $M$.  In the present paper we pay attention to the points
at infinity, i.e. the ``boundary'' $\partial M$ of the space $M$.  For each
point $e\in \partial M$ we introduce a new property of the action $\rho$,
analogous to $CC^{n-1}$ but defined using the filtration of $M$
by horoballs centered at $e$ rather than by the balls $B_r(a)$; we call
this ``$CC^{n-1}$ over $e$''.  Whereas in the previous situation the property
$CC^{n-1}$ was independent of the base point $a\in M$ it will now depend
in a delicate way upon the point $e\in \partial M$.  Therefore the subset 
of $\partial M$
\[
\Sigma^n(\rho) := \{e\in \partial M \mid\rho\ \mbox{ is }\ CC^{n-1}
\ \mbox{ over } e\},
\]
becomes interesting in its own right, and in this paper we investigate
its structure.

The isometric action $\rho$ of $G$ on $M$ induces a topological action of $G$ on
$\partial M$ when $\partial M$ is given the (compact) ``cone'' 
topology\footnote{The space $\partial M$ has two topologies:  the usual 
compact ``cone'' topology and the ``Tits distance'' topology which is usually
not compact but is complete and CAT$(1)$.  See \S11 for details.}, and
$\Sigma^n(\rho)$ is invariant with respect to this action.  However, it may
happen that the closure, $\cl_{\partial M}(Ge)$, of an orbit $Ge \subseteq
\Sigma^n(\rho)$ does not lie in $\Sigma^n(\rho)$.  Define
\[
\mathop\Sigma\limits^\circ{^n}(\rho) := \{e\in \partial M\mid 
\mbox{cl}_{\partial M}(Ge) \subseteq \Sigma^n(\rho)\}.
\]
There are important instances in which $\mathop\Sigma\limits^\circ{^n}(\rho)
= \Sigma^n(\rho)$ but when that is not the case
$\mathop\Sigma\limits^\circ{^n}(\rho)$ is easier to compute.  In more
detail:
\begin{itemize}
\item{}
If $\partial M$ is given the Tits distance topology then 
$\mathop\Sigma\limits^\circ{^n}(\rho)$ is an open subset\footnote{The set
$\mathop\Sigma\limits^\circ{^n}(\rho)$ is not, in general, open in the cone
topology; see the Remark following Theorem F in \S10.} of $\partial M$.  
In fact, $\{(\rho,e)\mid e\in \Sigma^n(\rho)\}$
is open in $\mbox{Hom}(G,\mbox{Isom}(M)) \x \partial M$. See Theorem F.

\item{} The set $\{\rho\mid \Sigma^n(\rho) = \partial M\}$ is open in 
Hom$(G$, Isom$(M))$.  See Theorem G.  Note that, for such $\rho$, 
$\Sigma^n(\rho) = \mathop\Sigma\limits^\circ{^n}(\rho)$.  

\item{} If $M$ is almost geodesically complete then $\Sigma^n(\rho) =
\partial M$ if and only if $\rho$ is uniformly $CC^{n-1}$ in the
sense of [BG$_{\mbox{I}}$].  See Theorem H. 

\item{} There is a description of 
$\mathop\Sigma\limits^\circ{^n}(\rho)$ in terms of a dynamical condition
on any free contractible $G$-CW-complex $X$ over $M$; see Theorem E.  This is
the key to all our results on $\Sigma^n(\rho)$ and can be used to compute
$\mathop\Sigma\limits^\circ{^n}(\rho)$ in specific cases.
\end{itemize}

When $M$ is Euclidean space $\EE^m$ then $\partial M$ is a sphere $S^{m-1}$
and the two topologies on $\partial M$ coincide; if the $G$-action $\rho$ on 
$M = \EE^m$ is by translations then the induced action on $\partial M$ is
trivial.  In this situation our results were previously known:  the two 
sets $\mathop\Sigma\limits^\circ{^n}(\rho)$ and $\Sigma^n(\rho)$ coincide
and have a description in terms of the {\em homotopical geometric invariant}
$\Sigma^n(G)$ of the group $G$, on which a considerable literature
exists.  More specifically, if we take $M$ to be the vector space $G/G'
\otimes \RR$ endowed with an inner product, and $\rho$ 
to be the action of $G$ on $M$
given by left translation then $\Sigma^n(\rho)$ has a direct 
interpretation as $\Sigma^n(G)$.  See Theorem I.  $\Sigma^n(G)$ has been
computed in many cases:  see \S\S10.6 and 10.7(A). 
Explicit computation of $\Sigma^n(\rho)$ in
non-Euclidean examples is still in its 
infancy, but we will be able to comment on the natural action of $G =
\SL_2(\ZZ_S)$, $S$ a natural number, on the hyperbolic plane (\S10.7(B)), 
and we
compute $\mathop\Sigma\limits^\circ{^n}(\rho)$ for all actions $\rho$ of
$G$ on locally finite simplicial trees (\S10.7(C)).  

The ideas in this 
two-part paper have evolved over the past twenty years.  $\Sigma^1(G)$
was first introduced in the special context of metabelian groups $G$
in [BS 80] as a tool to characterize finite presentability\footnote{Precisely:
a finitely generated metabelian group $G$ is finitely presented if and only if
for every point $e$ on the sphere at infinity either $e\in \Sigma^1(G)$
or $-e\in \Sigma^1(G)$.}  of $G$.  For arbitrary groups $G$, $\Sigma^1(G)$
was introduced in [BNS 87] and subsequently $\Sigma^n(G)$, $n\geq 2$,  
appeared in [BR 88] and [Re 88, 89].  Although these ``geometric invariants''
were defined for arbitrary groups $G$ of type $F_n$ they only referred
to translation actions on Euclidean spaces, and so $G$ was tacitly assumed
to have infinite Abelianization $G/G'$.  The present work emerged when we
tried to extend the scope of the theory to genuinely non-commutative
actions.  Not only the theorems but also the rather involved
techniques concerning finitary maps and their sheaves have
precursors in these earlier stages.

The numbering of sections and footnotes in this paper continues
that in [BG$_{\mbox{I}}$]. 

\setcounter{section}{9}

\newpage

\section{Outline, Main Results and Examples}

10.1 {\bf The boundary of a CAT(0)-space.} Let $(M,d)$ be a proper CAT(0)
space. A {\em geodesic ray} in $M$ is an isometric embedding $\gamma: [0,
\infty) \to M$. Two geodesic rays $\gamma,\gamma'$ are {\em asymptotic}
 if there is a
constant $r \in \RR$ such that $d(\gamma(t), \gamma'(t)) \leq r$
 for all $t$. The
set of all geodesic rays asymptotic to $\gamma$ is called the {\em endpoint
  of} $\gamma$ and denoted $ e = \gamma(\infty)$. The collection of all endpoints of
geodesic rays form the {\em boundary} 
$\partial M$ of $M$. Since $d$ is proper it is complete, so (see
[BrHa] or [Ho 97]) for every
 pair $(a,e) \in M \times \partial M$ there is a unique
geodesic ray $\gamma: [0,\infty) \to M$ with $\gamma(0) = a$ and $\gamma
(\infty) = e$. Hence there is a natural bijection between $\partial M$
 and the set of all geodesic rays emanating from a base
point, and $\partial M$ acquires the compact-open topology
of the latter via this bijection. 
This topology, which is independent of the choice of $a$, is compact and
metrizable.  It is called the {\em cone topology}.
 An action on $M$ by isometries induces a topological
action on $\partial M$.  $\partial M$ is assumed to carry the cone topology
except when another interesting topology, the Tits distance topology (see
\S11.4), is explicitly mentioned. 

Associated to each geodesic ray $ \gamma$ of $M$ is its {\em Busemann function}
$\beta_\gamma: M \to \RR$ (see \S 11.2) and for each $s \in \RR$ the associated
{\em horoball} $HB_s(\gamma):= \beta^{-1}_\gamma([s,\infty))$. Horoballs
 ``centered''
at $\gamma$ play a role analogous to that of balls centered at $a \in
M$. Indeed, $HB_s(\gamma) = {\rm cl}_M(\bigcup \{B_{t-s}(\gamma(t))|s < t\})$. 
Horoballs are contractible.
~\\

10.2 {\bf $\bf CC^{n-1}$ over end points.} In parallel with \S2.2, 
$M$ is a proper CAT(0) space,
$G$ is a group, $X^n$ an $n$-dimensional $(n-1)$-connected
free $G$-CW-complex such that $G\backslash X^n$ is finite, $\rho: G \to
\mbox{Isom}(M)$ is an isometric action and $h: X \to M$ is a control
function (i.e. a $G$-map). We pick a
geodesic ray $\gamma: [0,\infty) \to M$ and its endpoint $e =
 \gamma(\infty)$, and
we write $X_{(\gamma,s)}$ for the largest subcomplex of $X$ lying in
$h^{-1}(HB_s(\gamma))$. We say that $X$ {\em is controlled} $(n-1)$-{\em
  connected} $(CC^{n-1})$ {\em in the direction} $\gamma$ (with respect to
$\rho$) if for any horoball $HB_s(\gamma)$ and $-1\leq p\leq n-1$ there
exists $\lambda \geq 0$ such that every map $f : S^p \to X_{(\gamma,s)}$
extends to a map $\tilde f : B^{p+1} \to X_{(\gamma,s-\lambda)}$.  The
number $\lambda$ depends on the horoball $HB_s(\gamma)$ and is called a {\em
lag}.  In what follows, the lag will often turn out to be constant (i.e.
independent of $HB_s(\gamma)$; dependent only on $e$), and then we will
speak of a {\em constant lag}.  When $p = -1$ this
says that each $X_{(\gamma,s)}$ is non-empty.

The
property ``$X$ is $CC^{n-1}$ in the direction $\gamma$'' is shown in \S12.1
to be a property of the endpoint $e$ rather than
the ray $\gamma$. Moreover we have an Invariance Theorem (Theorem 12.1) as in
[BG$_{\mbox{I}}$]
which shows that this property is also independent of the choice
of $X$ and of $h$,
i.e., it is a property of the action $\rho$. So, if $X$ is $CC^{n-1}$ in
the direction $\gamma$ we will say that $\rho$ is $CC^{n-1}$ {\em over} (or 
{\em in the direction}) $e = \gamma(\infty)$.

So far the two conditions $CC^{n-1}$
over points $a \in M$, defined in \S2.2, and endpoints $e \in \partial M$ look
strictly analogous. But there is a striking difference in their
behaviour with respect to
changing $a$ or $e$: if $\rho$ is $CC^{n-1}$ over some $a \in M$ it is
also $CC^{n-1}$ over any other point of $M$. Not so for the property
$CC^{n-1}$ over $e \in \partial M$. Hence it becomes
interesting to study the subset of $\partial M$
\[ 
\Sigma^n(\rho) = \{e\; |\; \rho\; \mbox{is}\; CC^{n-1}\; \mbox{over}\; e\}.
\]

Following the tradition of [BS 80], [BNS 87], [BR 88] and [Re 88, 89] 
we call this the $n^{\mbox{th}}$ 
{\em (homotopical) geometric invariant of the action} $\rho$.
~\\

10.3 {\bf The dynamical subset.} From the definition it is clear that if
$e \in \partial M $ is in the geometric invariant $\Sigma^n(\rho)$, so
is the whole $G$-orbit $Ge$. The subset of $\Sigma^n(\rho)$
\[
\mathop\Sigma\limits^\circ{^n}(\rho) 
= \{e \in \partial M|\cl_{\partial M}(Ge)
\subseteq \Sigma^n(\rho)\} 
\]
deserves special attention. We call it the {\em dynamical subset} of
$\Sigma^n(\rho)$ because the following theorem characterizes it
in terms of a dynamical
condition in the free contractible $G$-CW-complex $X$ over $M$
(with cocompact $n$-skeleton) of \S 10.2. In the spirit of \S 2.6 we call a
cellular map $\phi: X^n \to X^n$ a {\em contraction towards} $e \in
\partial M$ if there is a number $\varepsilon > 0$ with
\[
\beta_\gamma h\phi(x) \geq \beta_\gamma h(x) + \varepsilon, \; \mbox{for
  all}\; x \in X^n,
\]
where $\gamma: [0,\infty) \to M$ is a geodesic ray with
$\gamma(\infty) = e$. 
We have the following
endpoint analogue of Theorem D: 
~\\

{\bf Theorem E.} {\em If} $\rho: G \to \Isom(M)$ {\em is an action of a
  group G of type} $F_n$ {\em on a proper} CAT(0) {\em space $M$ by isometries
  then the dynamical subset} $\mathop\Sigma\limits^\circ{^n}(\rho)$ {\em of}
$\Sigma^n(\rho)$ {\em can be characterized as}
\[
\mathop\Sigma\limits^\circ{^n}(\rho) = \{e|X^n\; \mbox{{\em admits a 
$G$-finitary contraction towards}}\; e\}.
\]

Note that 
$\mathop\Sigma\limits^\circ{^n}(\rho) = \partial M$ if and only
if $\Sigma^n(\rho) = \partial M$. 
The example of a Fuchsian group acting on the hyperbolic
plane (see \S 10.7(B)) shows that $\Sigma^n(\rho)$ and
$\mathop\Sigma\limits^\circ{^n}(\rho)$ are in general different but they
contain the same closed $G$-invariant subsets of $\partial M$.  
~\\

10.4  {\bf Openness results.}  The ``angular distance'' between two points
$e$ and $e'$ of $\partial M$ is the supremum over points $a\in M$ of the
angle between the geodesic rays $\gamma$ and $\gamma'$ representing $e$ and
$e'$ where $\gamma(0) = \gamma'(0) = a$.  This is a metric on $\partial M$
and the corresponding length metric is called the ``Tits distance'', denoted
$Td(e,e')$; see \S11.4 for full definitions.  In general, the topology on
$\partial M$ given by $Td$ is finer than the cone topology:  i.e. Id:
$(\partial M,Td) \to (\partial M$, cone topology) is continuous.  The
space $(\partial M,Td)$ is a complete CAT(1) metric space [BrHa; III 3.17].
Two extremes are represented by $M = \EE^k$ where the two topologies agree
giving $S^{k-1}$, and $M = \HH^k$ where $(\partial M,Td)$ is discrete while
$(\partial M$, cone topology) is $S^{k-1}$.  

We have openness theorems involving both topologies on $\partial M$.  The
easier of these involves Tits distance.  Recall from \S7.1 that the function
spaces Isom$(M)$ and Hom$(G$,Isom$(M))$ are given the compact-open topology.
With this understanding, the condition ``$e\in
\mathop\Sigma\limits^\circ{^n}(\rho)$'' is open 
with respect to both $\rho$ and $e$:
~\\

{\bf Theorem F.}  {\em $\{(\rho,e)\mid e \in
\mathop\Sigma\limits^\circ{^n}(\rho)\}$ is open in {\rm Hom}$(G,\mbox{\rm 
Isom}(M)) \x
\partial M$ when $\partial M$ carries the Tits distance topology. In
particular, for fixed $\rho$, $\mathop\Sigma\limits^\circ{^n}(\rho)$ is open
in $(\partial M,Td)$}.
~\\

{\bf Remark:} In general, neither $\Sigma^n(\rho)$ nor
$\mathop\Sigma\limits^\circ{^n}(\rho)$ is open in $(\partial M$,cone
topology). If $h$ is a non-trivial element of the free group of rank 2 and
$\rho$ is the action of $G := \langle h\rangle$ on the Cayley tree $T$ by
covering translations, then $\mathop\Sigma\limits^\circ{^n}(\rho)$ consists
of the two endpoints of the translation axis of $h$; this is not open in
$\partial T$.  For an example where $\Sigma^n(\rho)$ is not open see
\S10.7(B):  that example also shows that  
$\mathop\Sigma\limits^\circ{^n}(\rho)$ is not always the interior of
$\Sigma^n(\rho)$ in $(\partial M,Td)$.  We do not know if
$\Sigma^n(\rho)$ is always open in $(\partial M,Td)$.  In the
``classical'' case of [BS 80], [BNS 87], [BR 88] and 
[Re 88, 89] (discussed in \S10.7(A))
$\Sigma^n(\rho) = \mathop\Sigma\limits^\circ{^n}(\rho)$ and is open in
$S^{k-1}$.

For a different and deeper openness theorem we need $\partial M$ compact,
i.e. the cone topology:
~\\

{\bf Theorem G.}  {\em If $E$ is a closed subset of $(\partial M$, cone
topology) then $\{\rho \mid E\subseteq \Sigma^n(\rho)\}$ is an open subset of
{\rm Hom}$(G$,{\rm Isom}$(M,E))$.  
In particular, $\{\rho\mid \Sigma^n(\rho) = \partial
M\}$ is open in {\rm Hom}$(G$,{\rm Isom}$(M))$}.
~\\

Here, Isom$(M,E)$ is the space of isometries of $M$ which leave $E$
invariant.
~\\

10.5 {\bf Endpoints versus points in $M$.} 
We come to a theorem which relates the $CC^{n-1}$ property over all the
endpoints $e\in \partial M$ to the $CC^{n-1}$ property over a point $a\in M$.
This is the link between the present paper and [BG$_{\mbox{I}}$]; it requires
a mild assumption on $M$. The
(proper) CAT(0) space $(M,d)$ is {\em geodesically complete} if every
geodesic segment $[0,t] \to M$ can be extended to a geodesic ray
$[0,\infty) \to M$. We say the CAT(0) space $(M,d)$ is {\em almost
geodesically complete} if there is a number $\mu \geq 0$ such that for
any two points $a,b \in M$ there is a geodesic ray $\gamma$ with
$\gamma(0) = a$ such that $\gamma([0,\infty))$ meets $B_\mu(b)$. A
recent Theorem of P.~Ontaneda shows that this property is often
guaranteed in cases of interest.
~\\

{\bf Ontaneda's Theorem.} [On] {\em Let M be a non-compact proper}
 CAT(0)-{\em space such that} Isom$(M)$ {\em acts cocompactly. 
A sufficient
  condition for almost geodesic completeness is that the cohomology with
  compact supports} $H_c^*(M)$ {\em be non-trivial. This condition is
  satisfied whenever some subgroup of} Isom$(M)$ {\em acts cocompactly with
discrete orbits}\footnote{Recently, D. Farley has shown that this 
condition is also
satisfied if $M$ (as above) is an $M_\kappa$-complex with finite shapes (see
[BrHa] for the relevant definitions).}. 
~\\

{\bf Theorem H.} {\em If} $\rho: G \to \Isom(M)$ {\em is an isometric
  action on a proper and almost geodesically complete} CAT(0)-{\em space M
  then the following are equivalent}
\begin{enumerate}
\item[(i)] $\Sigma^n(\rho) = \partial M$
\item[(ii)] $\rho$ {\em is uniformly} $CC^{n-1}$ {\em over} $M$.
\end{enumerate}
~\\

We do not know whether the assumption that $M$ be almost geodesically
complete is necessary in Theorem H. 
It is not needed in the case $n = 0$, where
the Theorem says that $\Sigma^0(\rho) = \partial M$ if and only if
$\rho$ is cocompact.
It is an open problem as to whether there exists a non-compact CAT(0)
space $M$ which is not almost geodesically complete but admits a cocompact
group action by isometries.\footnote{The conclusion of Theorem B (``cocompact
and $CC^{n-1}$'' is an open condition) follows from
Theorems G and H.  This alternative proof of Theorem B has the 
disadvantage that it requires the (perhaps redundant) hypothesis that $M$
be almost geodesically complete.  It has the advantage of being technically 
simpler, in that $\partial M$ (with the cone topology)
plays the role previously played by the sphere of event
radius $R$ (\S5.3):  the compactness of this sphere is replaced by the
compactness of $\partial M$.  The reader who carefully compares these two 
proofs of Theorem B will see that this alternative proof is easier, avoiding
most of the content of \S\S6.3, 7.3 and 7.4.}  

The proofs of Theorems E, F,
G and H are completed in \S15.  
~\\

10.6 {\bf Fixed points and the BNSR-geometric invariant.} 
At this point we should relate our theory to the previous literature.
Let $S(G)$ be the sphere of non-zero (additive) characters on $G$ modulo 
positive
scalar multiplication:  when $\chi \in \Hom(G,{\RR})$ is a non-zero
character we denote by $[\chi]$ the point of $S(G)$ represented by $\chi$.
Identify $\RR$ with the group of translations of the Euclidean line and
interpret $\chi$ as a translation action of $G$ on $\RR$.  The preferred 
endpoint of $\RR$ is $\infty$.  Define
\[
\Sigma^n(G) := \{[\chi] \in S(G)\mid \chi\ {\rm is }\ CC^{n-1}\ {\rm  over }\ 
\infty\}.
\]
This is the {\it homotopical geometric invariant} of [BR 88] and [Re 88], 
and coincides with the invariant $\Sigma_{G'}$ of [BNS 87] when $n = 1$. 

Now we consider a CAT(0) space $M$ with a specified end point $e \in
\partial M$.  We denote by $\Isom(M,e)$ the group of isometries which fix
$e$.  If $f: M \to M$ is an isometry fixing $e$ then $f$ permutes
the horoballs centered at $e$, and so $f(H B_{(\gamma,r)}) = H
B_{(\gamma,r')}$ for $\gamma$ a geodesic ray with $\gamma(\infty) = e$
and $r \in \RR$. It is not difficult to observe that the map $r \mapsto
r'$ is a translation of $\RR$ (see Proposition 11.3) so that we have a
canonical map
\[
\kappa: \Isom(M,e) \to \RR = \Isom(\RR,\infty).
\]
Thus every $G$-action $\rho: G \to \Isom(M,e)$ induces a translation
action $\kappa\rho$ on $\RR$ -- in other words, an additive homomorphism
$\chi_{\rho,e}: G \to \RR$. Moreover, the Busemann function $\beta_\gamma: M
\to \RR$ is compatible with the actions $\rho$ and $\kappa\rho$ on $M$
and $\RR$, so if $h: X \to M$ is the control function chosen in \S 2.2
with respect to $\rho$ we can choose the composition $h_e = \beta_\gamma
h: X \to \RR$ as the control function for $\kappa\rho$. It is then immediate
from the definitions in \S 10.2 that $\rho$ is $CC^{n-1}$ over $e$ if
and only if $\kappa\rho$ is $CC^{n-1}$ over $\infty$. This proves 
~\\

{\bf Theorem I.} {\em Let} $\rho: G \to \Isom(M,e)$ {\em be an action
  fixing the point} $e \in \partial M$. {\em Then} $e \in
\Sigma^n(\rho)$ {\em (or, equivalently}, $e \in
\mathop\Sigma\limits^\circ{^n}(\rho)${\em ) if and only if}
$[\kappa\rho] \in \Sigma^n(G)$.

\hfill$\Box$
~\\

A more intimate relationship between $\Sigma^n(\rho)$ and $\Sigma^n(G)$
can be observed if $M = \EE^k$ is Euclidean $k$-space and $\rho: G \to
\mbox{Transl}(\EE^k)$ is an action on $\EE^k$ by translations. Then $\partial 
\EE^k$
is a $(k-1)$-sphere and is pointwise fixed by the $G$-action, so that $e
\mapsto [\chi_{\rho,e}]$ defines a map\footnote{$S(G) \cup \{0\}$ is a
quotient space of Hom$(G,\RR)$, $\{0\}$ is dense and $\mu$ is continuous.}
\[
\mu: \partial \EE^k = S^{k-1} \to S(G) \cup \{0\}.
\]
Let $N$ be the linear span of the orbit of 0, and let $N'$ be its orthogonal
complement in $\EE^k$. Then $\partial \EE^k = \partial N * \partial N'$
is the join of great subspheres in $\partial \EE^k$.  
The translation action $\rho$
of $G$ on $\EE^k$ 
restricts faithfully to a translation action $\rho_N$ on $N$, and
$\Sigma^n(\rho) = \Sigma^n(\rho_N) * \partial N' - \emptyset * \partial N'$.
The map $\mu$ respects the join operation, is injective on $\partial N$ and
maps $\partial N'$ to 0, so it maps $\Sigma^n(\rho_N)$ bijectively onto the
intersection of $\Sigma^n(G)$ with the great subsphere $\{[\chi]\mid \chi({\rm
ker}\ \rho) = 0\}$ of $S(G)$.  Thus the invariant $\Sigma^n(\rho)$ is determined
by $\Sigma^n(G)$.  In the special case where $\rho$ is cocompact and $G$ acts
properly discontinuously, $\mu$ maps $\Sigma^n(\rho)$ onto $\Sigma^n(G)$
homeomorphically.
~\\

10.7 {\bf Examples.}
~\\

A) {\bf Actions by Euclidean translations.} Let $G$ be a group of type
$F_n, n \geq 1$, and $M = G/G' \bigotimes \RR$ the real $G$-vector space
endowed with a Euclidean metric. The action $\rho: G \to
\mbox{Transl}(M)$ is induced by the left regular action of $G$ on 
$G/G'$. $\partial M$ is a $(k-1)$-sphere ($k =
\dim_{\RR}(G/G' \bigotimes \RR))$ which is pointwise fixed by
$\rho$. Thus, by \S 10.6, we have
\[
\mathop\Sigma\limits^\circ{^n}(\rho) = \Sigma^n(\rho) =
\Sigma^n(G) \subseteq S(G).
\]
The invariant $\Sigma^n(G)$ has been computed in many specific
cases.  We outline two situations which demonstrate the subtle behaviour of
$\Sigma^n(G)$.
~\\

A1) {\bf Right angled Coxeter groups.} Let $\Gamma$ be a graph. The group
$G(\Gamma)$ is given by the presentation whose generators are the
vertices of $\Gamma$ subject to the relations that adjacent vertices
commute. Every vertex $v \in \ver(\Gamma)$ gives rise to the
{\em coordinate hemispheres} $\pm H_v \subset S(G(\Gamma))$ where
\[
H_v:= \{[\chi] | \chi \in \Hom(G(\Gamma),\RR), \chi(v) > 0\}.
\]

Bestvina and Brady [BeBr 97] have shown that the diagonal character $\chi:
G(\Gamma) \to \RR$, defined by $\chi(v) = 1$ for all $v \in
\ver(\Gamma)$, lies in $\Sigma^n(G)$ if and only if the flag
complex\footnote{The flag complex of a graph is the simplicial complex
  whose $n$-simplices are the $(n+1)$-element sets of pairwise adjacent
  vertices} of $\Gamma$ is
$(n-1)$-connected. Meier, Meinert, and VanWyk [MMV 98] have extended this
to a complete description of $\Sigma^n(G(\Gamma))$; it is  
a finite union of finite intersections of coordinate hemispheres. 
For a simple geometric approach to
these results see [BuGo 99].
~\\

A2) {\bf Metabelian groups.} A group $G$ is said to be {\em metabelian}
if its commutator subgroup $G'$ is Abelian. For finitely generated
metabelian groups $G$ the invariant $\Sigma^1(G)$ has a handy description in
terms of valuations on fields [BS 81] which can be used for calculations.
In [BGr 84] it was used to prove that   
$\Sigma^1(G)$ is a rational polyhedral subset of $S(G)$, i.e., a finite
union of finite intersections of open hemispheres with rational
coordinate ratios for the spherical centres. It has been
conjectured that in the
metabelian case $\Sigma^n(G)$ is determined by $\Sigma^1(G)$; more
precisely, writing $\Sigma^n(G)^c$ for the complement of $\Sigma^n(G)$ in
$S(G)$, the $\Sigma^n$-{\em Conjecture} reads\footnote{
The $\Sigma^n$-Conjecture for metabelian groups $G$ came up in 1988.
It was based on the older $FP_n$-{\em Conjecture} (A finitely generated
metabelian group $G$ is of type $FP_n$
if and only if every $n$-point subset of $\Sigma^1(G)^c$ is contained in
an open hemisphere; see [B 81].).  Indeed, one can observe that if $G
\rightarrowtail \tilde G \twoheadrightarrow \ZZ$ 
is a short exact sequence then the $FP_n$-Conjecture
for $\tilde G$ implies the $\Sigma^n$-Conjecture for the rational points
in $\Sigma^n(G)$ up to antipodality.}: $G$ being of type $F_n$,
\[
\begin{array}{lcl}
\Sigma^n(G)^c &= &\{[\chi] \in S(G) \mid \chi = \chi_1 + \cdots +
\chi_k\\
& &\mbox{ with } k\leq n \ \mbox{ and each }\ [\chi_i] \in \Sigma^1(G)^c\}
\end{array}
\] 
This conjecture has been
verified when the commutator subgroup $G'$ is virtually
torsion-free of finite rank by H.~Meinert [Me 96] and when $G'$ is torsion
and of Krull dimension 1 as a $\ZZ G/G'$-module
by D.~Kochloukova [Ko 96]. Moreover, Kochloukova
has also proved it for general split\footnote{i.e. the commutator
  subgroup $G' \unlhd G$ has a complement in $G$.} metabelian groups and
$n = 2$.  We will see examples in C) below.
~\\

B) {\bf $\bf\SL_2(\RR)$-action on the hyperbolic plane.} 
Let $G$ be a discrete subgroup of $SL_2(\RR)$ (a Fuchsian group) and let
$\rho : G \to \Isom(\HH^2)$ be the natural action of $G$ on the
hyperbolic plane $\HH^2$.  If this action has 
a fundamental domain of finite area, then (see [BG 98]) the complement of
$\Sigma^0(\rho)$ in $\partial\HH^2$ is precisely the set of all
parabolic fixed points of $G$; or, equivalently, the orbits of the 
points $e \in
\partial\HH^2$ which lie on the boundary of a fundamental polygon. It is
not difficult to observe that, in fact, $\Sigma^n(\rho) =
\Sigma^0(\rho)$ for all $n$. Since the action of $G$ on $\partial\HH^2$
has only dense orbits we have
$\mathop\Sigma\limits^\circ{^n}(\rho) = \emptyset$ for all $n$.

A special case of a Fuchsian group is $\SL_2(\ZZ)$. If we use the upper
half plane model for $\HH^2$, the boundary is $\partial\HH^2 = 
\RR \cup\{\infty\}$ and $\Sigma^0(\rho)^c$ consists of the single orbit 
$\SL_2(\ZZ)\infty = \QQ \cup \{\infty\} =: \overline{\QQ}$. 

Let $S =
\{p_1, \ldots, p_s\}$ be a finite set of $s$ (different) prime numbers in $\NN$
and let $\ZZ_S$ be the subring of $\QQ$ generated by $(p_1p_2 \ldots
p_s)^{-1}$. We let $\rho: \SL_2(\ZZ_S) \to \Isom(\HH^2)$ be the
natural action of $\SL_2(\ZZ_S)$ on the upper half plane by Moebius
transformations. This action has dense orbits so, by Theorem H,  
$\mathop\Sigma\limits^\circ{^n}(\rho)$ is empty or all of
$\partial\HH^2$, and $\rho$ is $CC^{n-1}$ over $\infty \in \partial\HH^2$ 
if and only if $\rho$ is uniformly $CC^{n-1}$ in the sense of \S3. 
For this $\rho$ we propose: 
~\\

{\bf Conjecture.}
\[
\Sigma^n(\rho)^c = \left\{ \begin{array}{lll}
\emptyset & \mbox{if} & n < s\\
\overline{\QQ} & \mbox{if} & n \geq s\end{array} \right.
\]
 
{\em In particular}
\[
\mathop\Sigma\limits^\circ{^n}(\rho)^c = \left\{ \begin{array}{lll}
\emptyset & \mbox{if} & n \leq s\\
\partial\HH^2 & \mbox{if} & n > s \end{array} \right.
\]
~\\

The conjecture has been verified in special cases, in particular for the
case $s = 1$. Part of the interest of this conjecture has been discussed in
Footnote 11, but it is now given in a sharper form.  Eventually one would
hope for a similar statement for $SL_2(O_S)$ where $O_S$ is the ring of
$S$-integers in an algebraic number field.  
~\\

C)  {\bf Tree actions}.  Let $T$ be an infinite
locally finite tree and let $\rho : G \to
\Isom(T)$ be a cocompact action of $G$ by simplicial automorphisms.  Then, by
Bass-Serre theory, $G$ is the fundamental group of a finite graph of
groups $(\Gamma,{\cal G})$, where $\Gamma = G\backslash T$ and
${\cal G}$ is the system of edge and vertex stabilizers along a fundamental
transversal of $T$.  The edge stabilizers are of finite index in the
vertex stabilizers since $T$ is locally finite.  
Following [B 98] we define the {\em finiteness length} of $G$ [resp. ${\cal
G}$] to be $\fl G := \sup\{k\mid G$ is of type $F_k\}$ [resp. $\fl {\cal G}
:= \inf\{\fl H\mid H\in {\cal G}\}]$, and the {\em connectivity length} of a
character $\chi : G \to \RR$ to be $\cl(\chi) := \sup \{k\mid k\leq \fl G$
and $[\chi] \in \Sigma^k(G)\}$. In this case, $\fl {\cal G} = \fl H$ for any 
$H \in {\cal G}$.

We begin by noting three elementary facts.  First, $T$
is almost geodesically complete so we may apply Theorem H.  
Secondly, if the fixed point set $(\partial T)^G$ is a proper subset of
$\partial T$ then it is either empty or a singleton; for if there are two
singleton orbits then $\partial T$ consists of just those points.  Thirdly,
any orbit consisting of more than one point is dense, so that its closure is
$\partial T$.  It follows that if $S$ is the union of closures of
orbits---and $\mathop\Sigma\limits^\circ{^n}(\rho)$ is such a set---then $S =
\partial T$ or is a singleton or is empty.
It is well known that $\fl{\cal G} \leq \fl G$, so 
{\em if} $(\partial T)^G = \eset$ {\em then}
\begin{enumerate}
\item[(10.1)]
\[
\mathop\Sigma\limits^\circ{^n}(\rho) = \left\{ \begin{array}{lll} 
\partial T & \mbox{if} & 0 \leq n \leq \fl {\cal G}\\
\eset & \mbox{if} & \fl{\cal G} < n \leq \fl G \end{array} \right.
\]
\end{enumerate}

There remains the case when the $G$-tree $T$ has exactly one
fixed end $e$.
For such a tree we have the associated non-zero character
$\chi_{\rho,e} : G\to \RR$ of \S10.6 measuring the shift towards $e$. 
Thus if $n\leq \fl{\cal
G}$, $\Sigma^n(\rho)$ is defined, and by Theorems A and H, $\Sigma^n(\rho) =
\partial T$.  In particular, $e\in \Sigma^n(\rho)$ and therefore, by Theorem
I, $[\chi_{\rho,e}] \in \Sigma^n(G)$, implying $n\leq \cl(\chi_{\rho,e})$.
In summary:  $\fl{\cal G} \leq \cl(\chi_{\rho,e}) \leq \fl G$ and 
\begin{enumerate}
\item[(10.2)]
\[
\mathop\Sigma\limits^\circ{^n}(\rho) = \left\{ \begin{array}{lll} 
\partial T & \mbox{if} & 0 \leq n \leq \fl {\cal G}\\
\{e\} & \mbox{if} & \fl{\cal G} < n \leq \cl(\chi_{\rho,e})\\
\eset & \mbox{if} & \cl(\chi_{\rho,e}) < n \leq \fl G \end{array} \right.
\]
\end{enumerate}

By a {\em non-trivial rooted} $G$-{\em tree} we mean a (simplicial) $G$-tree
with no fixed point and with a unique fixed end $e$ (neither cocompactness nor
local finiteness are assumed here).  Every such $G$-tree $\rho : G \to
\mbox{Isom}(T)$ has an associated non-zero
character $\chi_{\rho,e} : G \to \RR$, defining an associated point
$[\chi_{\rho,e}] \in S(G)$. When $G$ is finitely generated then for
any non-trivial locally finite rooted $G$-tree there is a
(cocompact) Bass-Serre tree of an ascending HNN extension of $G$, i.e.,
\begin{enumerate}
\item[(10.3)]
\[
G = \langle B,t\mid tbt^{-1} = \theta(b),b\in B\rangle, 
\]
\end{enumerate}
where $\theta : B\rightarrowtail B$ is an injective endomorphism and
$\theta(B)$ has finite index $\geq 2$ in $B$, so that the two trees 
have the same associated
point of $S(G)$; i.e. one associated character is a positive
multiple of the other\footnote{More precisely, there is a
$G$-invariant ``superdivision'' $T_1$ of a $G$-invariant subtree $T'$ of $T$
(i.e. $T'$ is a subdivision of $T_1$) such that $G\backslash T_1$ consists
of one vertex and one edge.  This $T' : = \cup\{gA(h)\mid g\in G\}$ where
$A(h)$ is the translation axis of a hyperbolic element $h\in G$. There is
always a hyperbolic element under these hypotheses.}. 

We will examine the parameters $\fl {\cal G}$, $\cl(\chi_{\rho,e})$ and
$\fl G$ occurring in (10.2) when $G$ belongs to a special class of 
finitely generated groups.
For this, we define (for any character $\chi : G \to \RR$ on any finitely
generated group)
\[
m(\chi) := \sup\{k\mid\chi\neq \chi_1 + \cdots + \chi_k\ \mbox { with each
}\ [\chi_i] \in \Sigma^1(G)^c - \{[\chi]\}\}.
\]
Our special class is the class of finitely generated
{\em MFPR groups}, i.e. metabelian groups
of finite Pr{\"u}fer rank.  Recall that this means:  the commutator subgroup
$G'$ is abelian with finite torsion and finite torsion free rank
(dim$_{\QQ}(G'\otimes \QQ) < \infty$).  For such groups there is no
difference between type $F_n$ and type $FP_n$ ([BS 80, Theorem 5.4]).  Both the
$\Sigma^n$-Conjecture and the $FP_n$-Conjecture (see \S10.6) are known to
hold for MFPR groups ([\AA  86], [Me 96]).  Denoting the zero character by $0_G$
and an arbitrary character by $\chi$, these theorems can be restated in
terms of the function $m$:  
\begin{enumerate}
\item[(10.4)]
\[
\left\{ 
\begin{array}{ll}
\fl G &= m(0_G)\\
\cl(\chi) &= \min\{m(\chi),m(0_G)\}.
\end{array}
\right.
\]
\end{enumerate}
If $G$ splits as in (10.3), the character $\chi_{\rho,e}$ of the
corresponding tree decomposition is given by $\chi_{\rho,e}(B) = 0$ and
$\chi_{\rho,e}(t) = -1$. According to [Me 96; Theorem 2.6], 
{\em provided} $B$
{\em is finitely generated}, 
\begin{enumerate}
\item[(10.5)]
\[
\fl B = \min\{m(\chi_{\rho,e}),m(-\chi_{\rho,e}),m(0_G)\}.
\]
\end{enumerate}

The requirement that $B$ be finitely generated is not a serious restriction
if the MFPR group $G$ is finitely presented. In that case, for any 
non-trivial rooted
$G$-tree corresponding to a decomposition (10.3) there is another such,
having the same associated point of $S(G)$, in
which $B$ is finitely generated [BS 78; Theorem A]\footnote{In fact there is a 
map of $G$-trees from the latter onto the former.}.
Summarizing:  in this nice case ($G$ a finitely presented MFPR group and $B$
finitely generated) (10.2) can be rewritten:
\begin{enumerate}
\item[(10.6)]
\[
\mathop\Sigma\limits^\circ{^n}(\rho) = \left\{
\begin{array}{ll} 
\partial T &\mbox{if }\ 0\leq n \leq
\min\{m(\chi_{\rho,e}),m(-\chi_{\rho,e}),m(0_G)\}\\
\{e\} &\mbox{if }\ \min\{m(\chi_{\rho,e}),m(-\chi_{\rho,e}),m(0_G)\} < n \\
&\quad \leq \min\{m(\chi_{\rho,e}),m(0_G)\}\\
\ \ \eset &\mbox{if }\  \min\{m(\chi_{\rho,e},m(0_G)\} < n \leq m(0_G).
\end{array} \right.
\]
\end{enumerate}
Since the definition of $m(\chi)$ involves geometric features of the whole
set $\Sigma^1(G)^c$, it is only a slight exaggeration to conclude from
(10.6) that one must know $\Sigma^1(G)^c$ in order to know the
$\mathop\Sigma\limits^\circ{^n}$-properties of $\rho$.  

What can $\Sigma^1(G)^c$ be?  If $G$ is a finitely generated MFPR group,
then a precise answer can be given:  $\Sigma^1(G)^c$ is a finite set of
rational points in $S(G)$; moreover, given $k$ and any finite rational
subset $A$ of $S^k$ there is a finitely generated MFPR group $G$ with $S(G)
= S^k$ and $\Sigma^1(G)^c = A$; and in order to ensure that $G$ is
finitely presented, one simply chooses $A$ so that it does not contain
diametrically opposite points.  See [BS 81; Example 2.6].  

At this point we recall a theorem of K.S. Brown [Br 87$_{\mbox{II}}$] which
implies\footnote{Brown's theorem asserts that $\Sigma^1(G)^c$ consists of
all $[-\chi_{\rho,e}]$ such that $\rho$ is a $G$-action on an $\RR$-tree
with no fixed point and a unique fixed end.  For the rational points of
$S(G)$, simplicial (but not necessarily
locally finite) trees suffice.} that for any
finitely generated group $G$ the image of the map $\rho \mapsto
[\chi_{\rho,e}]$ from the set of all non-trivial rooted $G$-trees into $S(G)$
consists of all the rational points of $\Sigma^1(G)^c$.  It follows from our
discussion that if we only consider $G$-trees as in (10.3) with $B$ finitely
generated (call these ``special'') then, provided $G$ is finitely presented,
(i) there are enough special $G$-trees to map onto the rational points of
$\Sigma^1(G)^c$; (ii) those mapping to the same point have the same
$\mathop\Sigma\limits^\circ{^n}$-properties, and (iii) the
$\mathop\Sigma\limits^\circ{^n}$-properties of the various special $G$-trees
are interdependent in a manner dictated by the location of the finite set
$\Sigma^1(G)^c$ in $S(G)$.   

A final remark:  Among the choices of finitely generated MFPR groups $G$
and $[\chi] \in
\Sigma^1(G)^c$ we always have $m(\chi) = m(0_G)$ or $m(0_G) - 1$ (and both
can occur), but we
can make $m(\chi)$ and $m(\chi) - m(-\chi)$ as large as we like by 
choosing $k$ and $A$,
above, suitably. Thus, in (10.1) we can
achieve $\mathop\Sigma\limits^\circ{^n}(\rho) = \eset$ only when $n = \fl G
< \infty$, while we can achieve $\mathop\Sigma\limits^\circ{^n}(\rho) =
\partial T$ and $\mathop\Sigma\limits^\circ{^n}(\rho) = \{e\}$ in
arbitrarily large ranges of $n$.  We do not know if
$\mathop\Sigma\limits^\circ{^n}(\rho) = \eset$ can occur in a larger range
when $G$ is not an MFPR group.

\newpage

\section{Further Technicalities on CAT(0) spaces}

11.1 {\bf More on $\bf \partial M$.} In order to relate the topology
of the CAT(0) space $M$ to that of its boundary
$\partial M$ we use the following notation\footnote{For proofs of the
assertions about CAT(0) spaces in this section see [BrHa].}: 
A {\em generalized geodesic
  ray} $\gamma: [0,\infty) \to M$ is either a geodesic ray in the sense
of \S 10.1 or a geodesic segment $[0,\mu] \to M$ extended by
$\gamma([\mu,\infty)) = \gamma(\mu)$; it {\em starts} at $\gamma(0)$ and {\em
  ends} at $\gamma(\infty):= \gamma(\mu)$. 

Let $R$ be the function space consisting of all generalized geodesic
rays $\gamma: [0,\infty) \to M$ with the compact open topology. (i.e.,
uniform convergence on compact sets). Let $\hat{M}$ denote the disjoint
union $M \cup \partial M$, and $\varepsilon: R \dpo \hat{M}$ the
endpoint map $\varepsilon(\gamma):= \gamma(\infty)$. For each $a \in M$
the map $\varepsilon$ has a continuous canonical section $\sigma_a: \hat{M} \pfm
R$, where $\sigma_a(e)$ is the unique generalized geodesic ray starting
at $a \in M$ and ending at $e \in \hat{M}$. This shows that $\hat{M}$
with the quotient topology induced by $\varepsilon$ is homeomorphic to
the subspace $R_a:= \sigma_a(\hat{M}) = \{\gamma \in R|\gamma(0) = a\}$ of
$R$. Restricting $\sigma_a$ to the subspace $M \subset \hat{M}$ shows
that the topology of $M$ inherited from $\hat{M}$ is the topology given
by the metric on $M$; and $\partial M$ is homeomorphic to the subspace
$\partial R_a:= \sigma_a(\partial M) \subset R$ consisting of all
geodesic rays emanating from $a$. 
~\\

{\bf Proposition 11.1.} $\hat{M}$ {\em is a compact metrizable space
  containing M as a subspace,} $\partial M = \hat{M}-M$ {\em is a
  nowhere dense subset and for every open set} $U \subset \hat{M}$, {\em
  the inclusion} $U-\partial M \pfm U$ {\em is a homotopy equivalence.}

\hfill$\Box$
~\\

When a
group $G$ acts on $M$ by isometries there is an obvious extension to an
action of $G$ on $\hat{M}$ by homeomorphisms: If $e \in \partial M$ is
the endpoint of the geodesic ray $\gamma: [0,\infty) \to M$, and $g \in
G$, then $g\gamma$ is a geodesic ray and we put $ge:=
(g\gamma)(\infty)$. The stabilizer of $e \in \partial M$ in $G$ is
denoted $G_e$.
~\\

11.2 {\bf Review of Busemann functions and horoballs.} 
Associated to a generalized
geodesic ray $\gamma: [0,\infty) \to M$ is its Busemann function
$\beta_\gamma: M \to \RR$ defined by
\[
\beta_\gamma(b) = \lim\limits_{t\to\infty}(d(\gamma(0),\gamma(t)) -
d(b,\gamma(t)), \qquad b \in M.
\]
In the case when $\gamma(\infty) \in M$ and $\mu$ is the smallest
non-negative number with $\gamma(\mu) = \gamma(\infty)$ we find
$\beta_\gamma(b) = \mu-d(b,\gamma(\mu))$. In the case when
$\gamma(\infty) \in \partial M$ one shows, by the triangle inequality,
that the map $t \mapsto t - d(b,\gamma(t))$ is monotone increasing and
bounded above by $d(\gamma(0),b)$, so $\beta_\gamma(b)$ is well
defined. 
~\\

{\bf Proposition 11.2.} (a) {\em If} $\gamma$ {\em and} $\gamma'$ {\em
  are generalized geodesic rays with the same endpoint then}
$\beta_\gamma(\cdot)-\beta_{\gamma'}(\cdot)$ {\em is constant. In
  particular the difference $\beta_\gamma(a)-\beta_\gamma(b)$ {\em
    {\em depends only on} $a,b$, {\em and} $e = \gamma(\infty)$.\\
(b) {\em The function} $\beta: R \times M \to \RR,\; \beta(\gamma,b) =
\beta_\gamma(b)$ {\em is continuous.}

\hfill $\Box$
~\\

If $\gamma$ is a generalized ray and $s \in \RR$, then
the (closed) {\em horoball} $HB_s(\gamma)$ is
$\beta_\gamma^{-1}([r,\infty))$. When $s \geq 0$ this has a more
geometric interpretation, namely $HB_s(\gamma) =
\cl_M(\cup\{B_{t-s}(\gamma(t)) | s < t\})$. If the ray $\gamma:
[0,\infty) \to M$ is degenerate with $e = \gamma(\infty) \in M$ and
$d(\gamma(0),e) = \mu$ then $HB_s(\gamma)$ is precisely $B_{\mu-s}(e)$.
~\\

11.3 {\bf $\bf G$-actions and Busemann functions.} 
Let $G$ be a group acting on $M$ by
isometries, and let $\gamma: [0,\infty) \to M$ be a 
geodesic ray. For $g \in G$ and $a \in M$ we observe that
$\beta_\gamma(ga) - \beta_\gamma(a)$ depends only on $g,a$ and $e =
\gamma(\infty)$. So putting $\psi_e(g,a):=
\beta_\gamma(ga)-\beta_\gamma(a)$ defines a map $\psi_e: G \times M \to
\RR$. This map satisfies the equation
\[
(11.1)\quad \psi_e(gh,a) = \psi_e(g,ha) + \psi_e(h,a), \qquad g,h
\in G.
\]
To see how $\psi_e(g,a)$ depends on $a$ we compute
\[
\begin{array}{lcl}
\psi_e(g,a)-\psi_e(g,b) & = &
\beta_\gamma(ga)-\beta_\gamma(a)-\beta_\gamma(gb)+\beta_\gamma(b)\\
& = &
(\beta_{g^{-1}\gamma}(a)-\beta_\gamma(a))-(\beta_{g^{-1}\gamma}
(b)-\beta_\gamma(b)),
\end{array}
\]
where we have used the obvious fact that $\beta_{gx}(ga) =
\beta_x(a)$. Now suppose $g$ lies in the stabilizer $G_e$ of
$e$. Then $g^{-1}\gamma$ and $\gamma$ have the same endpoint and
since $\beta_{g^{-1}\gamma}(\cdot)-\beta_\gamma(\cdot)$ is constant we
find $\psi_e(g,a) = \psi_e(g,b)$ for all $g \in G_e,\; a,b \in M$. By
(11.1) we get:
~\\

{\bf Proposition 11.3.} {\em The function} $\chi_e: G_e \to \RR$ {\em
  defined by} $\chi_e(g) = \beta_\gamma(ga)-\beta_\gamma(a)$ {\em is 
 an additive  homomorphism. In other words,} 
$g \cdot r:= r+\chi_e(g)$ {\em defines a
  translation action of} $G_e$ {\em on} $\RR$ {\em which makes the
Busemann function} $\beta_\gamma: M \to \RR$ a $G_e$-{\em map}.

\hfill $\Box$

11.4 {\bf The Tits distance.}  We need some preliminary definitions.  If
$\Delta$ is a geodesic triangle in $M$ with vertices $a$, $b$ and $c$, and
if $\Delta'$ is the comparison triangle in the plane, the {\em angle at} $a$
in $\Delta$ is the Euclidean angle at $a$ in $\Delta'$; compare \S2.1.  The
{\em angle} between geodesic rays $\gamma$ and $\gamma'$ starting at the
same point $a\in M$ is $\angle_a(\gamma,\gamma') := \ds{\lim_{\epsilon \to 0}} \
\ds{\sup_{t,t' < \epsilon}}\{$ the angle at $\gamma(0)$ with vertices
$\gamma(t)$ and $\gamma'(t')\}$.  The {\em angular distance} between points
$e$ and $e'$ of $\partial M$ is $\angle(e,e') := \ds{\sup_{a \in
M}}\{\angle_a(\gamma,\gamma')\}$ where $\gamma$ and $\gamma'$ (as above)
represent $e$ and $e'$ respectively.  This is a metric called the {\em
angular metric} on $\partial M$: it is bounded above by $\pi$.  The {\em
Tits distance}, $Td(e,e')$ between $e,e' \in M$ is the inf of lengths of
rectifiable (with respect to the angular metric) paths\footnote{If there are
no such paths the Tits distance is defined to be $\infty$; this extension of
the usual notion of ``metric'' causes no problems.} in $\partial M$.
These two metrics are complete and are locally isometric to one another.
Hence they define the same topology on $\partial M$.  The identity map
$(\partial M,Td) \to (\partial M$, cone topology) is continuous.  The Tits
distance gives rich structure to $\partial M$; in particular $(\partial
M,Td)$ is a complete CAT(1) space [BrHa; III 3.17] where ``CAT(1)'' is defined
in the same way as CAT(0) but using comparison triangles in $S^2$ rather
than in the plane.  

\newpage
\section{$CC^{n-1}$ over Endpoints}

12.1 {\bf Invariance Theorem.} As in \S 3.1 $M$ is a proper CAT(0)
space, $X$ is a free left $G$-complex with $G\backslash X^n$
  finite and $\rho: G \to \Isom(M)$ is a left action of $G$ on $M$ by
  isometries. By Proposition 3.1 we can choose a $G$-equivariant
  control function $h: X \to M$.

A {\em decreasing filtration} of $X$ is a collection of subcomplexes 
$(K_t)_{t \in \RR}$ of $X$ such that $K_s \supseteq K_t$ whenever $s \leq t$.
Filtrations $(K_t)$
and $(L_t)$ will be called {\em equivalent} if there is a constant
$\mu \geq 0$ such that $K_{t-\mu} \supseteq L_t \supseteq K_{t+\mu}$ for
all $t$.

The geodesic ray $\gamma$ defines a filtration of $X$ by the
subcomplexes $\{X_{(\gamma,s)}|s \in \RR\}$, where $X_{(\gamma,s)}$ is
the largest subcomplex of $X$ lying in $h^{-1}(HB_s(\gamma))$. The
equivalence class of this filtration only depends on
$\gamma(\infty) = e \in \partial M$ since, by Proposition 11.2,
$X_{(\gamma',t)} = X_{(\gamma,t+{\rm constant})}$ when $\gamma'$ is
asymptotic to $\gamma$. And it does not depend on the $G$-map $h$ since,
referring to Proposition 3.1, for any $A \subset M\; h_2^{-1}(A)$ lies
in $h_1^{-1}$(the $t$-neighbourhood of $A$), where $t =
\sup\{d(h_1(x),h_2(x))|x \in X^n\}$.  Moreover, 
the property that $X$ be $CC^{n-1}$ over $e$ is,
in fact, also independent of the choice of $X$, so that it is really a
property of the action $\rho$ and the point
$e \in \partial M$. This
is covered by the following more general Invariance Theorem
~\\

{\bf Theorem 12.1.} {\em Let G be of type} $F_n$, {\em let Y be a
  cocompact $n$-dimensional $(n-1)$-connected rigid $G$}-CW-{\em complex 
such that the
  stabilizer of each p-cell is of type} $F_{n-p},\; p \leq n-1$. {\em
  Let h:} $Y \to M$ {\em be a G-map and let} $e \in \partial M$. {\em
  The property that Y be} $CC^{n-1}$ {\em over e is independent of the
  choices of $\gamma$, of Y and of the control function h.}
~\\

The proof of Theorem 12.1 carries over, mutatis mutandis, from the proof
of Theorem 3.3.

In \S10.2 we remarked that the lag $\lambda$ depends on the horoball
$HB_s(\gamma)$.  By Proposition 11.2, if $\gamma'$ is asymptotic to $\gamma$
then $HB_s(\gamma) = HB_{s'}(\gamma')$ for some $s'$.  We mean that
$\lambda$ is not a function of $\gamma$ and $s$ but rather depends only on
the set $HB_s(\gamma)$. 
~\\

12.2 {\bf $\bf CC^{-1}$ in all directions.} Here we prove
~\\

{\bf Theorem 12.2} {\em The action} $\rho: G \to \Isom(M)$ {\em is}
$CC^{-1}$ {\em in all directions} $e \in \partial M$ {\em if and only
  if} $\rho$ {\em is cocompact.}
~\\

{\em Proof.} We may take $n = 0$ and consider a finite subset $F \subset
X^0$ with $GF = X^0$. If $\rho$ is cocompact then there is a compact set
$C \subset M$ such that $GC = M$ and $h(F) \cap C \not=
\emptyset$. Since $C$ has finite diameter, each horoball must contain a
translate of $C$, hence some point of $h(X^0) = Gh(F)$; i.e. $X$ is
$(-1)$-connected in all directions.

It remains to show that if $\rho$ is not cocompact then there is a
horoball $HB_s(\gamma)$ disjoint from $h(X^0)$. There is a sequence $(b_m)$
in $M$ such that $\inf d(h(X^0), b_m) \geq m$ for all $m$.  Let $\epsilon
> 0$.  There are points $a_m \in h(X^0)$ such that $\inf d(h(X^0),b_m) \leq
d(a_m,b_m) \leq \inf d(h(X^0),b_m)+\epsilon$,  
and translating each $b_m$ with
the $G$-action if necessary we may assume that $a_m \in h(F)$. Pick $a
\in h(F)$. Since the metric $d$ is proper the space $R_a$ of all
generalized geodesic rays emanating from the base point $a$ is compact. 
The points $b_m \in M$ are
represented in $R_a$ by the unique generalized geodesic rays $\omega_m:
[0,\infty) \to M$ starting in $a$ and ending at $\omega_m(\infty) = b_m$. The
sequence $(\omega_m)$ has a limit point $\gamma \in R_a$, and we may assume
it converges to $\gamma$. If $\gamma(\infty) \in M$ then the sequence
$(d(a,b_m))$ is bounded which contradicts the choice of $(b_m)$. Hence
$\gamma(\infty) \in \partial M$, i.e., $\gamma: [0,\infty) \to M$ is a
geodesic ray. As $(\omega_m)$ converges to $\gamma$ we find for each $k \in
\NN$ a number $N(k) \geq k$ with $d(\omega_{N(k)}(k), \; \gamma(k)) <
\epsilon$. We abbreviate $c_k:= \omega_{N(k)}(k)$ and observe that, 
as well as $d(\gamma(k), c_k) < \epsilon$, we have also $d(a,c_k) = k$, since
$c_k$ is the parameter-$k$-point on the geodesic segment from $a$ to
$b_{N(k)}$ which has length $d(a,b_{N(k)}) + \epsilon\geq d(a_{N(k)}, b_{N(k)})
\geq N(k) \geq k$.

We claim that for each point $c \in h(X^0)$ the inequality
\[
d(a,c_k) \leq d(a,a_{N(k)})+d(c,c_k) + \epsilon
\]
holds. Indeed, this follows from
\[
\begin{array}{lcl}
d(a,b_{N(k)}) & \leq & d(a,a_{N(k)}) + d(a_{N(k)},b_{N(k)})\\
& \leq & d(a,a_{N(k)}) + d(c,b_{N(k)}) +\epsilon\\
& \leq & d(a,a_{N(k)}) + d(c,c_k) + d(c_k,b_{N(k)}) +\epsilon
\end{array}
\]
by subtracting $d(c_k,b_{N(k)})$ on either side. It follows that for
each $c \in h(X^0)$ and each $k \in \NN$

\[
\begin{array}{lcl}
d(c,\gamma(k)) & \geq & d(c,c_k) - d(c_k,\gamma(k))\\
& \geq & d(a,c_k) - d(a,a_{N(k)}) - d(c_k,\gamma(k))-\epsilon\\
& \geq & k - \diam h(F) - 2\epsilon.
\end{array}
\]

This shows that none of the points of $h(X^0)$ is contained in the
horoball $HB_{{\rm diam}\, h(F)+2\epsilon}(\gamma)$.

\hfill $\Box$

\newpage

\section{Finitary Contractions Towards Endpoints}

Throughout this section $X$ is a contractible $G$-CW-complex with finite
stabilizers and cocompact $n$-skeleton, $M$ is a CAT(0) metric space, an action
of $G$ on $M$ by isometries is given, and
$h: X \to M$ is a $G$-equivariant control function. 
The main technical result in [BG$_{\mbox{I}}$] was devoted to 
characterizing ``$CC^{n-1}$ over $a \in M$'' in terms of the existence of
a finitary contraction towards $a$. We now consider finitary contractions
towards an endpoint $e \in \partial M$ in connection with the
$CC^{n-1}$ property over $e$.
~\\

13.1 {\bf Shift and contractions towards $\bf e \in \partial M$.}
 As in \S 5.2 we
consider cellular maps $f: D(f) \to X$, where $D(f)$ is a
subcomplex of $X$. The {\em shift of $f$ towards} $e \in \partial M$ (or
{\em in the direction} $e$) is defined to be the continuous function
$\sh_{f,e}: D(f) \to \RR$ given by
\[
\sh_{f,e}(x):= \beta_\gamma hf(x) - \beta_\gamma h(x), 
\]
where $\gamma$ is a geodesic ray representing $e$. By Proposition 11.2(a)
this is independent of the choice of $\gamma$. By the definition of
$\beta_\gamma$ and the triangle inequality the shift is bounded by the
displacement function 
\[
(13.1) \qquad |\sh_{f,e}(x) | \leq \alpha_f(x), 
\]
hence, just as for the shift towards $a \in M$, the shift function has a
global bound $\|f\|$ if $f$ is a bounded map\footnote{Recall from \S5.1 that
$\alpha_f(x) := d(h(x),hf(x))$ and $||f|| := \sup \alpha_f(D(f))$.}.

By the {\em guaranteed shift towards} $e \in \partial M$ we mean
\[
\gsh_e(f):= \inf \sh_{f,e}(D(f)).
\]
Just as in \S 5.4 we find
\[
\gsh_{ge}(gf) = \gsh_e(f), \qquad \mbox{all}\; g \in G.
\]

A cellular map $\phi: X \to X$ is said to be a {\em contraction towards} $e
\in \partial M$ (or {\em in the direction} $e$) if $\gsh_e(\phi) > 0$.
~\\

13.2 {\bf From contractions to $\bf CC^{n-1}$.} The endpoint versions
of the propositions in \S 5.3 have easier statements and proofs. 
Corresponding to
Proposition 5.3 we have the observation that
\[
(13.2)\qquad \gsh_e(\phi^m) \geq m\gsh_e(\phi)
\]
for each contraction $\phi: X \to X$ and each $m \in \NN$. Corresponding
to Proposition 5.5 we have
~\\

{\bf Proposition 13.1} {\em If} $X^n$ {\em
  admits a finitary contraction} $\phi: X^n \to X^n$ {\em towards} $e \in
\partial M$ {\em and} $\epsilon := {\rm gsh}_e(\phi)$  
 {\em then there is a number}
$\lambda \geq 0$ {\em and a cellular deformation} $\psi: X^n \times
  [0,\infty) \to X^{n+1}$ {\em satisfying the Lipschitz condition}
\[
\beta_\gamma h\psi(x,s_2)-\beta_\gamma h\psi(x,s_1) \geq
(s_2-s_1)\varepsilon-\lambda
\]
{\em whenever} $s_1 \leq s_2$, $x \in X^n$, {\em and} $\gamma$
{\em is a geodesic ray with} $\gamma(\infty) = e$.
~\\

{\em Proof.} As in the proof of Proposition 5.5 we define $\psi_0: X^n
\times I \to X^{n+1}$ to be a finitary homotopy $\Id_{X^n}\simeq \phi$
and take $\psi(x,t): = \psi_0(\phi^m(x),s)$, where $x \in X^n,\; t \in
[0,\infty)$, and $m$ is an integer with  $s = t-m \in I$. Define 
$\sh_{\psi_{0,e}}(x) := \ds{\sup_{t\in I}}\ \sh_{\psi_{0,e}(\cdot,t)}(x)$.
This time we
 note that for $y:= (\phi^m(x),t_1)$ and
$z:= (\phi^m(x),t_2)$ with $0 \leq t_1 \leq t_2 \leq 1$ and $x \in X^n$,
\[
\begin{array}{lcl}
|\beta_\gamma h\psi_0(y)-\beta_\gamma h\psi_0(z)| & \leq & |\beta_\gamma
h\psi_0(y) - \beta_\gamma
h(\phi^m(x))| + |\beta_\gamma h(\phi^m(x))-\beta_\gamma h\psi_0(z)|\\
& \leq & |\sh_{\psi_0,e}(\phi^m(x))| + | \sh_{\psi_0,e}(\phi^m(x))|\\
& \leq & 2 \|\psi_0\|.
\end{array}
\]
If $s_i:= m_i+t_i$ as above, with $s_1 \leq s_2$, it follows that
\[
\begin{array}{lcl}
\beta_\gamma h\psi(x,s_2)-\beta_\gamma h\psi(x,s_1) & = & \beta_\gamma
h\psi_0(\phi^{m_2}(x),t_2)-\beta_\gamma h\psi_0(\phi^{m_1}(x),t_1)\\
& \geq & \beta_\gamma h\psi_0(\phi^{m_2}(x),1)-\beta_\gamma
h\psi_0(\phi^{m_1}(x),0)-4\|\psi_0\|\\
& \ = & \beta_\gamma h\phi^{m_2+1}(x)-\beta_\gamma
h\phi^{m_1}(x)-4\|\psi_0\|\\
& = & \sh_{\phi^{m_2-m_1+1},e}(\phi^{m_1}(x))-4\|\psi_0\|\\
& \geq & \gsh_e(\phi^{m_2-m_1+1})-4\|\psi_0\|\\
& \geq & (m_2-m_1+1)\varepsilon-4\|\psi_0\|, \quad \mbox{by (13.2)}\\
& \geq & (s_2-s_1)\varepsilon-4\|\psi_0\|.
\end{array}
\]
So we can choose $\lambda$ to be $4\|\psi_0\|$.

\hfill $\Box$
~\\

By a straightforward adaptation of the proof of Theorem 5.6 we get
~\\

{\bf Theorem 13.2} {\em If X is contractible with finite 
  stabilizers and cocompact n-skeleton then the existence of a finitary
  contraction} $\phi: X^n \to X^n$ {\em towards} $e \in \partial M$ {\em
  implies that X is} $CC^{n-1}$ {\em over e with constant lag} $\lambda =
4\|\psi_0\|$, {\em where} $\psi_0$ {\em is any finitary homotopy} 
$\Id_{X^n} \simeq \phi$.

\hfill $\Box$
~\\

There is a subtle variation of Theorem 13.2 which did not arise in the
parallel case in \S 5 but will be needed in the inductive proof of the
main result (Proposition 14.3).
~\\

{\bf Corollary 13.3} {\em The existence of a finitary contraction} 
$\phi: X^n \to X^n$ {\em towards}
  $e \in \partial M$ {\em also implies that if X has cocompact}
  $(n+1)$-{\em skeleton and is} $CC^n$ {\em over e then it is so with
constant lag} $\lambda = 4\|\psi_0\|$.
~\\

{\em Proof.} Assume $X$ is $CC^n$ over $e$ with lag $\lambda_0$ and
consider a map $f: S^n \to X_{(e,t)}$. Put $\varepsilon: = \gsh_e
(\phi)$ and $r:= \varepsilon^{-1}\lambda_0$, and compose the deformation
$\psi$ of Proposition 13.1 with the map $(f \times \Id): S^n \times
[0,r] \to X_{(e,t)} \times [0,r]$. By the Lipschitz condition of $\psi$
this yields a homotopy $H: S^n \times [0,r] \to X_{(e,t-\lambda)}$
between $f$ and a map $f' = H(\cdot,r): S^n \to
X_{(e,t+\lambda_0-\lambda)}$. By the $CC^n$-property of $X$ over $e$,
with lag $\lambda_0,\; f'$ extends to a map $\tilde{f'}: B^{n+1} \to
X_{(e,t-\lambda)}$. Now, the annulus $S^n \times [0,r]$ and the ball
$B^{n+1}$ can be glued together along $S^n \times \{r\}$ to a
topological $(n+1)$-ball $\tilde{B}$, and the union $\tilde{f} = H \cup
\tilde{f'}$ is a map $\tilde{f}: \tilde{B} \to X_{(e,t-\lambda)}$
extending $f$.

\hfill $\Box$
~\\

13.3 {\bf Passing to the closure of $G$-orbits.} It is clear that if
$\phi: X^n \to X^n$ is a finitary contraction towards $e \in \partial M$
then $g\phi$ is a finitary contraction towards $ge$, so that, by
Theorem 13.2, $X$ is $CC^{n-1}$ over each $e'$ in the $G$-orbit of
$e$, with a uniform constant lag $\lambda$. We improve this by establishing
~\\

{\bf Proposition 13.4} {\em The existence of a finitary
  contraction} $\phi: X^n \to X^n$ {\em towards} $e \in \partial M$ {\em
  implies that X is} $CC^{n-1}$ {\em over every} $e'$ {\em in}
$\cl_{\partial M}(Ge)$, {\em the closure of the G-orbit of e. Moreover,
  if} $\psi_0$
{\em is a finitary homotopy} $\Id_{X^n}\simeq \phi$ {\em then any
  number} $\lambda' > 4 \|\psi_0\|$ {\em is a constant lag, uniform 
for all} $e' \in
\cl_{\partial M}(Ge)$.
~\\

In the proof of Proposition 13.4 we will need the following
~\\

{\bf Lemma 13.5} {\em Let M be a} CAT(0) {\em space and} $\varepsilon,r$
{\em positive numbers. Then any number} $R > r(1+2r/\varepsilon)$ {\em has
the property that whenever} $\gamma,\gamma': [0,\infty) \to M$ {\em are
geodesic rays starting at} $c \in M$ {\em then}
\[
|\beta_\gamma(p) - \beta_{\gamma'}(p) | < 2 \varepsilon + d(\gamma(R),
\gamma'(R))
\]
{\em for every} $p \in B_r(c)$.
~\\

{\em Proof.} Applying Lemma 6.4 with $a = \gamma(t), b = \gamma(R)$ and
$p' = c$ yields
\[
|d(c,\gamma(t)) - d(p,\gamma(t)) - d(c,\gamma(R)) + d(p,\gamma(R)) | <
\varepsilon.
\]
Hence passing to the limit $t \to \infty$ we find
\[
|\beta_\gamma(p) - R + d(p,\gamma(R)) | < \varepsilon,
\]
and similarly for $\gamma'$. Subtracting the two inequalities and using
the triangle inequality yields the lemma.

\hfill $\Box$
~\\

{\em Proof (of Proposition 13.4).}  We use the notation of the proof
of Proposition 13.1.  The proof of Theorem 13.2 starts with $0\leq p\leq n-1$,
a map $f : S^p \to X^p$, its extension to $f_1 : B^{p+1} \to X^{p+1}$, and 
a number $t$ such that $f(S^p) \subset X_{(\gamma,t)}$.  The deformation
$\psi$ is used to move $f_1$ to $f'_1 : B^{p+1} \to X^{p+1}$, where 
$f'_1 := \psi(f_1(\cdot),T)$, $T$ being sufficiently large that $f'_1(B^{p+1})
\subseteq X_{(\gamma,t-\lambda)}$.  The Lipschitz condition on $\psi$ shows
that if $\epsilon := {\rm gsh}_e(\phi)$ then $T =$ diam $hf_1(B^{p+1})
/\epsilon$ will do.  Let $r$ be the diameter over $M$ of the $[0,T]$ 
track of $f_1(B^{p+1})$, i.e. $r:=$ diam $h\psi(f_1(B^{p+1})\x [0,T])$.
The final map $\tilde f : B^{p+1} \to X^p$ satisfies diam $h\tilde f(B^{p+1})
\leq r$ and $\tilde f(B^{p+1}) \subset X_{(\gamma,t-\lambda)}$.  

The numbers $\lambda = 4||\psi_0||$, $\epsilon = \gsh_e(\phi)$, $T$ and $r$
remain unchanged when $e$ and $\phi$ are replaced by $G$-translates $ge$ and
$g\phi$.  This shows that for every $g\in G$ the given map $f : S^p \to
X^p$ can be extended to a map $\tilde f_g : B^{p+1} \to X^{p+1}$ such
that (i)~diam~$h\tilde f_g(B^{p+1}) \leq r$ and (ii) $f(S^p) \subset
X_{(g\gamma,s)}$ implies $\tilde f_g(B^{p+1}) \subset X_{(g\gamma,
s-\lambda)}$.  In words, (ii) says that the extension $\tilde f_g$ of $f$
has lag $\lambda$ with respect to $ge$.

Now choose a base point $c\in hf(S^p)$.  For all $g\in G$ we have $h\tilde f_g
(B^{p+1}) \subset B_r(c)$.  Let $e' \in \cl(Ge)$  
be represented by a geodesic ray $\gamma':
[0,\infty) \to M$ with $\gamma'(0) = c$. Let $\nu > 0$ be an arbitrary
small number, choose $R > r(1+2r/\nu)$, and choose $g \in G$ such that
if $ge$ is represented by the geodesic ray $\overline{\gamma}:
[0,\infty) \to M$ with $\overline{\gamma}(0) = c$, then $d(\gamma'(R),
\overline{\gamma}(R)) < \nu$. Lemma 13.5 shows that
\[
|\beta_{\gamma'}(u) - \beta_{\overline{\gamma}}(u) | < 3\nu,\;
\mbox{for all}\; u \in B_r(c).
\]
This shows that the extension 
$\tilde{f}_g: B^{p+1} \to X^{p+1}$ of $f:
S^p \to X^p$ has lag $\lambda + 3\nu$ with respect to $e'$.

\hfill $\Box$
~\\

13.4 {\bf Contractions in a set of directions ${\bf E \subseteq \partial
    M}$}. We conclude this section with a further improvement on
Proposition 13.4 and Corollary 13.3.
~\\

{\bf Proposition 13.6} {\em Let} ${\cal F}: X^n
\to X^n$ {\em be a locally finite G-sheaf and let} 
$E \subseteq \partial M$ {\em be a
closed set of endpoints with the property that for each} $e \in E$ {\em the
sheaf} ${\cal F}$ {\em has a cross section} $\phi_e: X^n \to X^n$ {\em
with} $\inf\{\gsh_e(\phi_e)|e \in E\} = \varepsilon > 0$. {\em Then
there is a number} $\lambda$ {\em such that the following hold:} \\
(a) X {\em is} $CC^{n-1}$ {\em over each} $e' \in \cl_{\partial M}(GE)$
 {\em with constant lag} $\lambda$\\
(b) {\em If X has cocompact} $(n+1)$-{\em skeleton and is} $CC^n$ {\em
  over} $e' \in GE$ {\em then it is so with lag} $\lambda$.
~\\

{\em Proof.} The lag for an individual endpoint $e \in E$ which we found
in Corollary 13.2 was $4\|\psi_{0,e}\|$ where $\psi_{0,e}$ was a finitary
homotopy $\Id_{X^n}\simeq \phi_e$. Such a homotopy is obtained as a
cross section of a homotopy of sheaves ${\cal H}: {\cal F}_0 \simeq {\cal
  F}$, where ${\cal F}_0$ stands for the sheaf consisting of all
identity maps $\Id_D$ with $D \subseteq X^n$ the domain of a member of
${\cal F}$; see Propositions 4.8 and 4.10.  Since $\psi_{0,e}$ is a cross
section of ${\cal H}$, we have $||\psi_{0,e}|| \leq ||{\cal H}|| :=
\sup\{||H||\mid H\in {\cal H}\}$.  Therefore
$4\|{\cal H}\|$ is a uniform lag for all $e \in E$, and Proposition
13.4 shows that any number $\lambda > 4\|{\cal H}\|$ is a uniform lag
for all $e' \in \cup\{\cl\ Ge\mid e \in E\}$.  Since $E$ is closed in
$\partial M$ this set is $\cl_{\partial M}GE$.   

\hfill $\Box$

\newpage

\section{From $CC^{n-1}$ over Endpoints to Contractions}

Let $X$, $M$ and $h$ be as in \S13.  We aim to construct
finitary contractions from $CC^{n-1}$ assumptions over endpoints.
~\\

14.1 {\bf Vertex shift and defect of sheaves.} 
Let ${\cal F}: X \leadsto X$ be a locally finite homotopically
closed $G$-sheaf. In \S 5.4 and \S 5.5 we introduced the maximal 
guaranteed vertex shift $\mu_a({\cal F}|\sigma)$ and the defect 
$d_a({\cal F}|\sigma)$ on a cell $\sigma$ towards a point $a \in M$. The same
definitions apply with $a \in M$ replaced by an endpoint $e \in \partial
M$.

Guaranteed vertex shift and defect in the direction $e \in \partial M$
are used to control the shift of a cross section $\phi: X \to X$ of
${\cal F}$ on $\sigma$ in the direction $e$. Just as in \S 5.5 one
introduces the total defect of ${\cal F}$ on $\sigma$ in the direction 
$e,\; \delta_e({\cal F}|\sigma)$, and one proves (compare Proposition 5.7
and Remark 5.8):
~\\

{\bf Proposition 14.1} {\em Let} $e \in \partial M$. {\em Then every
  locally finite homotopically
closed $G$-sheaf} ${\cal F}: X \leadsto X$ {\em has a cross
  section} $\phi: X \to X$ {\em with}
\[
\sh_{\phi,e}(x) \geq \mu_e({\cal F}|\sigma) - \delta_e({\cal F}|\sigma)
\]
{\em for each cell} $\sigma$ {\em of X and each} $x \in \sigma$; and $\phi$
{\em can be chosen to be a} $G_e$-{\em map}.

\hfill $\Box$
~\\

14.2. {\bf Controlled embedding of sheaves into homotopically 
closed sheaves.} We
shall need the following embedding result
~\\

{\bf Proposition 14.2} {\em Let E be a closed and G-invariant subset
  of} $\partial M$. {\em Assume X is} $CC^{n-1}$ {\em in all
  directions} $e \in E$ {\em with uniform constant lag} $\lambda \geq 0$ {\em
  depending only on E. Then for every} $\lambda' > \lambda$, {\em any
  locally finite G-sheaf} ${\cal F}: X^n \leadsto X^n$ {\em can be
    embedded in a homotopically closed locally finite G-sheaf}
  $\tilde{{\cal F}}$ {\em such that} $d_e(\tilde{\cal F}|\sigma) \leq
  \lambda'+\diam hC(\sigma)$, {\em for all cells} $\sigma$ {\em of}
  $X^n$, {\em and all} $e \in E$.
~\\

{\em Proof.} The argument is parallel to those in \S 6.2 and \S 6.4, in
fact considerably simpler since the technicalities of \S 6.3 do not
arise. Just as in \S 6.2 one starts by fixing an $n$-cell $\sigma$ of
$X$ and a cellular map $f: C(\mathop\sigma\limits^\bullet) \to X$ and
one proves that, given any $\varepsilon > 0$, a finite set $S(\sigma)$ of
cellular maps $C(\sigma) \to X$ extending $f$ can be chosen in such a
way that for each geodesic ray $\gamma$ with $\gamma(\infty) = e \in E$
there exist $\tilde{f}_e \in S(\sigma)$ with
\[
\gsh_\gamma(\tilde{f}_e) > \gsh_\gamma(f) - \diam hC(\sigma) - \lambda -
\varepsilon.
\]
Note that we have used the assumption that $E$ is compact and the fact
that the Busemann function $\beta_\gamma$ is continuous in $\gamma$
(Proposition 11.2).

We can now prove Proposition 14.2 by induction on $n$. When $n = 0$ the
condition on the defects $d_e({\cal F}|\sigma)$ is empty so we can refer
to Proposition 4.5. By induction, assume that ${\cal F}^{n-1}:= {\cal F}\mid
X^{n-1} : X^{n-1}
\leadsto X^{n-1}$ is homotopically closed with $d_e({\cal F}^{n-1}|\tau) \leq
\lambda'+\diam hC(\tau)$ for every cell $\tau$ of $X^{n-1}$ and all $e
\in E$. In order to embed ${\cal F}$ in a homotopically closed sheaf we adjoin
the sets $S(\sigma)$ constructed above to ${\cal F}$. The set $S(\sigma)$
depends on $\sigma$; but for the $G$-translates $g\sigma$ of $\sigma$ we
do not have to choose $S(g\sigma)$ anew but can put $S(g\sigma):=
gS(\sigma)$. This is because $\gsh_{ge}(gf) = \gsh_e(f)$; see \S
13.1. In this way we have embedded ${\cal F}$ in a sheaf $\tilde{{\cal
    F}}: X^n \leadsto X^n$ which is homotopically closed, locally finite and a
$G$-sheaf; and $\tilde{{\cal F}}$ does satisfy the required condition on
the defect.

\hfill $\Box$
~\\

14.3 {\bf Obtaining contractions towards end points.} 
Here is the converse to Proposition 13.6. 
~\\

{\bf Proposition 14.3} {\em If} $E \subseteq \partial M$ {\em is a
  closed $G$-invariant subset and X is} $CC^{n-1}$ {\em over}
{\em every} $e \in E$ {\em then there is a locally finite homotopically closed
  G-sheaf} ${\cal F}: X^n \leadsto X^n$ {\em which admits, for each} $e
\in E$, {\em a cross section} $\phi_e: X^n \to X^n$ {\em with}
$\inf\{\gsh_e(\phi_e)|e \in E\} > 0$.
~\\

{\em Proof.} First we consider the case $n = 0$. For each vertex $v$ of
a finite fundamental domain $F \subseteq X^n$ and each end point $e \in
E$ we use the $CC^{-1}$ hypothesis to pick a map $f^v_e: \{v\} \to X^0$ 
with $\gsh_e(f^v_e) = h
f^v_e(v)-h(v) > 0$. For $e' \in E$ sufficiently close to $e,\;
\gsh_{e'}(f^v_e) > 0$. Since $E$ is compact we find a finite sheaf
${\cal F}_0: F^0 \leadsto X^0$ with $\mu_e({\cal F}_0|v) > 0$ for all $e
\in E$ and all $v \in F^0$. Since $F^0$ and the stabilizers of the vertices
$v \in
F^0$ are finite there are only finitely many group elements $g \in G$ with
$gF^0 \cap F^0 \not= \emptyset$. By adding to ${\cal F}_0$
the maps $gf$ for all $f \in
{\cal F}_0$ and all $g \in G$ with $gD(f) \in F^0$ we find a finite
$G$-saturated sheaf ${\cal F}'_0: F \leadsto X^0$ (in the sense of \S
4.2).

By Propositions 4.1 and 4.3 ${\cal F}^0 := G{\cal F}'_0: X^0 \leadsto
X^0$ is a locally finite $G$-sheaf. Since it contains ${\cal F}'_0$ and
since $\mu_{ge}({\cal F}^0|gv) = \mu_e({\cal F}^0|v)$ we have
$\mu_e({\cal F}^0|x) > 0$ for all $e \in E$ and all $x \in X^0$ -- in
fact $\inf\mu_e({\cal F}^0|X^0) > 0$. The required cross section
$\phi_e$ of ${\cal F}^0$ exists by Proposition 14.1.

We proceed by induction on $n$, assuming that a homotopically
closed locally finite
$G$-sheaf ${\cal F}: X^{n-1} \leadsto X^{n-1}$ has already been constructed,
with cross sections $\phi_e: X^{n-1} \to X^{n-1}$ contracting
towards $e \in E$, and
$\varepsilon:= \inf\{\gsh_e(\phi_e)|e \in E\} > 0$. The assumptions of
Proposition 13.6(b) are then fulfilled (with $n$ replaced by $n-1$), so
there is a uniform constant lag for the
$CC^{n-1}$-assumptions on $X$ in all directions $e \in E$. In this
situation Proposition 14.2 applies and yields a number $D$ with the
property that every locally finite $G$-sheaf ${\cal F}': X^n \leadsto
X^n$ can be embedded into a locally finite homotopically closed
$G$-sheaf $\tilde{\cal F}': X^n \leadsto X^n$ 
with $d_e(\tilde{{\cal F}}'|\sigma) \leq D$ for all $e \in E$, and
all cells $\sigma$ of $X^n$.

We may alter the domain of ${\cal F}$ and write ${\cal F}: X^n \leadsto
X^n$, noting that ${\cal F}\mid X^{n-1} : X^{n-1} \leadsto X^{n-1}$ 
(but perhaps not ${\cal F}$) is homotopically 
closed. Let $m$ be a positive integer such that $m \cdot \varepsilon >
D$. The $m$-th iterate ${\cal F}^{(m)}: X^n \leadsto X^n$ is a locally
finite $G$-sheaf and thus, by the previous paragraph, can be embedded
into a locally finite homotopically closed $G$-sheaf $\tilde{{\cal F}}:
X^n \leadsto X^n$ with $d_e(\tilde{{\cal F}}|\sigma) \leq D$ for all $e
\in E$ and all cells $\sigma$ of $X^n$. Then the maps $\phi^m_e: X^{n-1}
\to X^{n-1}$ are cross sections of $\tilde{{\cal F}}$ and, by (13.2),
$\gsh_e\phi^m_e \geq m\gsh_e(\phi_e) \geq m \cdot \varepsilon$. From the
definition of the defect $d_e(\tilde{{\cal F}}|\sigma)$ it is clear that
each of the cross sections $\phi^m_e: X^{n-1} \to X^{n-1}$ can be
extended to a cross section $\tilde{\phi}_e: X^n \to X^n$ of the
homotopically closed sheaf $\tilde{{\cal F}}$ with
\[
\begin{array}{lcl}
\gsh_e(\tilde{\phi}_e|C(\sigma)) & \geq & \gsh_e(\phi^m_e|C(\sigma)) -
d_e(\tilde{{\cal F}}|\sigma)\\
& \geq & m\varepsilon - D > 0,
\end{array}
\]
for all cells $\sigma$ of $X^n$.

\hfill $\Box$
~\\

14.4 {\bf The main technical results.} 
~\\

{\bf Theorem 14.4} {\em If X is a contractible G-{\rm CW}-complex with finite
  cell stabilizers and cocompact n-skeleton and if $G$ acts on $M$ by
isometries then the following are
  equivalent for} $e \in \partial M$.
\begin{enumerate}
\item[(i)] $X^n$ {\em admits a G-finitary contraction towards e}
\item[(ii)] {\em X is} $CC^{n-1}$ {\em in all directions} $e' \in
  \cl_{\partial M}(Ge)$
\item[(iii)] {\em There is a locally finite homotopically closed
    G-sheaf} ${\cal F}: X^n \leadsto X^n$ {\em and a number}
  $\varepsilon > 0$ {\em with}
\[
 \mu_e({\cal F}|\sigma)-\delta_e({\cal F}|\sigma) \geq \varepsilon
\]
\end{enumerate}
{\em for all cells} $\sigma$ {\em of} $X^n$.
~\\

{\bf Addendum.} If (i)--(iii) {\em  hold then there is a constant lag} 
$\lambda \geq 0$ {\em uniform for all} $e' \in \cl(Ge)$ {\em in} (ii). 
~\\

{\em Proof.} (i) $\Rightarrow$ (ii) is covered by Proposition 13.4
which yields a uniform constant lag $\lambda$ as required in the Addendum. 
(ii) $\Rightarrow$ (i)
is covered by Proposition 14.3. The implication (iii) $\Rightarrow$ (i)
is clear from Proposition 14.1. It remains to show that (i) implies (iii).  We
know that (i) implies (ii) with uniform constant 
lag $\lambda$ as described in the
Addendum.  Let $\lambda' > \lambda$ and let $L > \lambda' + {\rm diam}
\ hC(\sigma)$ for all cells $\sigma$ of $X^n$.  Let $\phi : X^n \to X^n$ be a
finitary contraction towards $e$.  Then $\gsh_e(\phi) > 0$, and for some $m
\geq 1$ $\phi^m$ is a finitary contraction towards $e$ with $\gsh_e(\phi^m) >
L$; this uses Proposition 4.9.  Then ${\cal F}' := G\ {\rm Res}(\phi^m)$
is a locally finite $G$-sheaf with $\mu_e({\cal F}'\mid \sigma) > L$ for all
cells $\sigma$ of $X^n$.  By Proposition 14.2, using (ii) with uniform
constant lag
$\lambda$ on cl$_{\partial M}(Ge)$, ${\cal F}'$ can be embedded in a
homotopically closed locally finite $G$-sheaf ${\cal F}$ with $d_e({\cal F}\mid
\sigma) \leq \lambda' + {\rm diam}\ hC(\sigma) < L$ for all cells $\sigma$ of
$X^n$.  Since $\mu_e({\cal F}\mid \sigma) > L$ for all $\sigma$, (iii) holds.

\hfill $\Box$
~\\

{\bf Remark.}  For any cell $\sigma$ of $X$ and any $\cal F$ as above, the
numbers $\mu_e({\cal F}\mid \sigma)$ and $\delta_e({\cal F}\mid \sigma)$
depend continuously on $e\in \partial M$; this follows from Proposition
11.2; compare \S5.4.  Hence the inequality $\mu_e({\cal F}\mid \sigma) -
\delta_e({\cal F}\mid \sigma) > \frac{\epsilon}{2}$ which comes from (iii)
of Theorem 14.4 for given $e$ implies $\mu_{e'}({\cal F}\mid \sigma)
-\delta_{e'}({\cal F}\mid \sigma) > \frac{\epsilon}{2}$ for all $e$ in some
neighborhood of $e$.  However, since there are infinitely many $\sigma$ to
be considered one cannot conclude that
$\mathop\Sigma\limits^\circ{^n}(\rho)$ is open in $\partial M$ with respect
to the cone topology.  Counterexamples are found in \S\S10.4, 10.7(B) and
10.7(C).
~\\

{\bf Theorem 14.5} {\em If X is a contractible $G$-{\rm CW}-complex with finite
  cell stabilizers and cocompact n-skeleton and if $G$ acts on $M$ by
isometries then the following are
  equivalent for a G-invariant closed subset} $E \subseteq \partial M$.
\begin{enumerate}
\item[(i)] $X^n$ {\em admits G-finitary contractions} $\phi_e: X^n \to
  X^n$ {\em towards each} $e \in E$.
\item[(ii)] {\em X is} $CC^{n-1}$ {\em in all directions} $e \in E$.
\item[(iii)] {\em There is a locally finite homotopically closed
    G-sheaf} ${\cal F}: X^n \leadsto X^n$ {\em and a number}
  $\varepsilon > 0$ {\em with}
\[
\mu_e({\cal F}|\sigma) - \delta_e({\cal F}|\sigma) \geq \varepsilon
\]
{\em for all cells} $\sigma$ {\em of} $X^n$ {\em and all} $e \in E$.
\end{enumerate}
~\\

{\bf Addendum.} {\em If (i) -- (iii) hold then more is true.  There is a
locally finite homotopically closed $G$-sheaf ${\cal F} : X^n \rightsquigarrow
X^n$ such that the maps $\phi_e$ in (i) can be chosen to satisfy {\rm
Res}$(\phi_e) \subset {\cal F}$ for all $e$ and 
$\inf\{\gsh_e(\phi_e)|e \in E\}
> 0$; and (ii) holds with a uniform constant
lag $\lambda$ in all directions $e \in E$.} 
~\\
 
{\em Proof.} (i) $\Leftrightarrow$ (ii) follows from Theorem 14.4.  (iii)
$\Rightarrow$ (i) is immediate from Proposition 14.1. It remains to prove
that (i) and (ii) imply (iii).   Assuming (ii),
Proposition 14.3 gives the stronger form of (i) in the Addendum; and this
implies the stronger form of (ii) in the Addendum (uniform constant
lag $\lambda$) by
Proposition 13.6.  Thus in proving (i) $\Rightarrow$ (iii) we may assume the
stronger forms of (i) and (ii).  The proof runs parallel to the (i)
$\Rightarrow$ (iii) part of Proposition 14.4.  Let $L > \lambda' + {\rm diam}\
hC(\sigma)$ for all cells $\sigma$ of $X^n$. For each $e\in E$ let
$\phi_e : X^n \to X^n$ be a
finitary contraction towards $e$, each $\phi_e$ being a  
cross section of the sheaf ${\cal F}$ in the Addendum, and
$\inf\{{\rm gsh}_e(\phi_e)\mid e\in E\} > 0$.  By Propositions 4.7 and 4.9, each
iterate $\phi^m_e$ is a finitary contraction towards $e$ and is a cross section
of the $m$-fold composite sheaf ${\cal F}' := {\cal F}^{(m)}$. This is a
locally finite $G$-sheaf and, if $m$ is large enough, $\mu_e({\cal F}'\mid
\sigma) > L$ for all cells $\sigma$ of $X^n$ and all $e\in E$.  The rest of the
proof is similar to the corresponding proof in Theorem 14.4.  

\hfill $\Box$

\newpage

\section{Proofs of Theorems E-H}

\ \ \ \ \ 15.1  {\bf Dynamical characterization of}
$\mathop\Sigma\limits^\circ{^n}(\rho)$. 
~\\

{\em Proof of Theorem E}.  This is the equivalence of (i) and (ii) in
Theorem 14.4. 

\hfill$\Box$
~\\

15.2 {\bf Openness using the Tits distance topology}. 
Recall from \S5.1 that the norm of a cellular map $\phi : X^n \to X^n$ is
$||\phi|| := \sup\{d(h(x),h\phi(x)) \mid x\in X^n\}$, and that $\phi$ is
said to be bounded if $||\phi|| < \infty$.  Recall from \S10.3 that $\phi$
is a contraction towards $e\in \partial M$ if there is a number $\epsilon >
0$ such that $\sh_{\phi,e}(x) \geq \epsilon$ for all $x\in X^n$ (see \S13.1
for $\sh_{\phi,e}(x)$).  
~\\

{\bf Proposition 15.1}  {\em Let $\phi$ be a bounded contraction towards
$e \in \partial M$.  There is a neighborhood $N$ of $e$ in $(\partial M,Td)$
such that $\phi$ is also a bounded contraction towards every $e' \in N$}.  
~\\

{\em Proof}.  Let $\gamma$ and $\gamma'$ be geodesic rays representing $e$
and $e' \in \partial M$ such that $\gamma(0) = \gamma'(0)$.  Let $\alpha(t)$
be the angle at $\gamma(0)$ in the geodesic triangle whose other vertices
are $\gamma(t)$ and $\gamma'(t)$.  With notation as in \S11.4, it is shown
in [BrHa, III 3.4] that 
$\alpha(t) \leq \angle(e,e')$ and $\ds{\lim_{t\to\infty}}
\alpha(t) = \angle(e,e')$.  By Euclidean geometry $\sin \frac{\alpha(t)}{2} =
\frac{d(\gamma(t),\gamma'(t))}{2t}$, so 
\[
d(\alpha(t),\gamma'(t)) \leq 2t \sin \left(\frac{\angle(e,e')}{2}\right).
\]
The right hand side of this inequality is independent of the base point
$\gamma(0)$.

Let $\epsilon > 0$ be such that $\sh_{\phi,e}(x) \geq \epsilon$ for all $x\in
X^n$.  We apply Lemma 13.5 (with $\frac{\epsilon}{6}$ replacing $\epsilon$
in that lemma, and $r = ||\phi||$).  Letting $R = ||\phi||(1 +
\frac{12||\phi||}{\epsilon}) + \epsilon$, that lemma implies that for any
$x\in X^n$ and any $p \in B_{||\phi||}(h(x))$ we have:
\[
|\beta_\gamma(p) - \beta_{\gamma'}(p)| < \frac{\epsilon}{3} +
d(\gamma(R),\gamma'(R))
\]
where $\gamma(0) =\gamma'(0) = h(x)$. Thus if $e'$ is chosen so that
$\angle(e,e') < 2 \arcsin (\frac{\epsilon}{6R})$ we find 
\[
\begin{array}{l l}
|\sh_{\phi,e}(x) - \sh_{\phi,e'}(x)| &= |\beta_\gamma h\phi(x) -
\beta_{\gamma'}h\phi(x)|\\
 &\leq \frac{\epsilon}{3} + \frac{\epsilon}{3R} R = \frac{2\epsilon}{3}.
\end{array}
\]
This shows that $\sh_{\phi,e'}(x) > \frac{\epsilon}{3}$ for all $x\in X^n$.

\ \hfill $\Box$
~\\
 
Parallel to Proposition 7.3 we have (when $h_\rho$ is chosen continuously in
$\rho$)
~\\

{\bf Proposition 15.2} {\em Let} ${\cal F}: X \leadsto X$ {\em be a
  locally finite homotopically closed G-sheaf. Then} $\mu^\rho_e({\cal
  F}|\sigma)$ and $\delta^\rho_e({\cal F}|\sigma)$ {\em are jointly
  continuous in the variables} $(\rho,e) \in \Hom(G,\Isom(M)) \times
\partial M$, {\em where $\partial M$ has the cone topology}.
~\\

Of course this remains true if $\partial M$ has the finer Tits distance
topology. 
~\\

{\em Proof of Theorem F}.  For fixed $\rho$, openness in $e$ follows
from Proposition 15.1 and Theorem 14.4.  If $\rho$ is permitted to vary and
$h_\rho$ is chosen continuously in $\rho$, the same proof works, in view of
Proposition 15.2.  

\ \hfill $\Box$
~\\

15.3 {\bf Openness in $\rho$ using the cone topology on $\partial M$}.
~\\

{\em Proof of Theorem G}. 
Let $\rho \in \mbox{Hom}(G,\mbox{Isom}(M,E))$ be such that one and hence all
of the Conditions (i)-(iii) of Theorem 14.7 hold.  Since $\mu_{ge}({\cal
F}\mid g\sigma) = \mu_e({\cal F}\mid\sigma)$ and $\delta_{ge}({\cal F}\mid
g\sigma) = \delta_e({\cal F}\mid \sigma)$ for all $g\in G$ the inequality in
Condition (iii) need only be considered for a (finite) set of
representatives of the $G$-orbits of cells in $X^n$.  The parameter range
$e\in E$ is compact and hence the condition (with $\epsilon$ replaced by
$\frac{\epsilon}{2})$ remains true if $\rho$ is subject to a small
perturbation, in view of Proposition 15.2.  

\ \hfill $\Box$
~\\

15.4 {\bf Endpoints versus points in M.}  Theorem H follows from the
Invariance Theorem 12.1 and 
~\\

{\bf Theorem 15.3.}  {\em Let $X$ be a contractible $G$-{\rm CW}-complex with
finite cell stabilizers and cocompact $n$-skeleton.  If $h : X \to M$ is an
equivariant control function over an almost geodesically complete}
CAT(0)-{\em space $M$ then the following are equivalent}
\begin{enumerate}
\item[(i)]  $X$ {\em is $CC^{n-1}$ in all directions} $e\in \partial M$
\item[(ii)]  $X$ {\em is uniformly} $CC^{n-1}$ {\em over} $M$.
\end{enumerate}
~\\

{\em Proof.}  (ii) $\Rightarrow$ (i) is easy and is left to the reader.  We
prove (i) $\Rightarrow$ (ii).  Let ${\cal F}: X^n \leadsto X^n$ be the sheaf
given by (iii) of Theorem 14.5 for $E = \partial M$.  Let $r := ||{\cal
F}||$, let $\nu$ be an arbitrarily small positive number and choose $S > r(1
+ 2r/\nu)$.  Let $\mu \geq 0$ be the number given by the definition of
``almost geodesically complete'' in \S10.5.  We may assume $\nu \leq 1$ and
$\mu \geq 1$ and put $R := \mu S/\nu$.  Now we fix a base point $a\in M$ and
consider any $x\in X^n$ with $d(a,h(x)) \geq R$.  By assumption there is a
geodesic ray $\gamma : [0,\infty) \to M$ with $\gamma(0) = h(x)$ and
$d(a,\gamma(T)) \leq \mu$ for some $T\in [0,\infty)$.  We put $e :=
\gamma(\infty)$ and we compare the shifts $\sh_{f,a}(x)$ and $\sh_{f,e}(x)$
for all $f\in {\cal F}$.  

Since $R \geq S > r(1+2r/\nu)$ Lemma 6.4 applies.  We apply it, in fact,
twice, each time taking $c = p = h(x)$ and $p' = hf(x)$. The first
application is for the point $\gamma(t) \not\in B_S(c)$ and yields
\[
|\sh_{f,\gamma(t)}(x) - \sh_{f,\gamma(S)}(x)| < \nu.
\]
Taking the limit $t\to\infty$ this becomes
\[
(15.1)\quad |\sh_{f,e}(x) - \sh_{f,\gamma(S)}(x)| < \nu.
\]
The second application is for the point $a\not\in B_S(c)$ and yields
\[
(15.2) \quad |sh_{f,a}(x) - \sh_{f,b}(x)| < \nu,
\]
where $b$ is the point on the geodesic segment from $c$ to $a$ with $d(c,b)
= S$.  In order to get a bound on $d(\gamma(S),b)$ we have to compare the
two isosceles triangles $(a,c,\gamma(d(a,c)))$ and $(b,c,\gamma(S))$ with
the corresponding Euclidean ones.  Since $d(a,\gamma(d(a,c))) \leq
d(a,\gamma(T)) + |T-d(a,c)| \leq 2d(\alpha,\gamma(T)) \leq 2\mu$ we find
$d(\gamma(S),b)/S \leq 2\mu/d(a,c) \leq 2\mu/R$, and hence $d(\gamma(S),b)
\leq 2\nu$.  By (5.7) it follows that 
\[
(15.3) \quad |\sh_{f,\gamma(S)}(x) - \sh_{f,b}(x)| \leq 4\nu,
\]
and hence, from the conjunction of (15.1) - (15.3),
\[
(15.4) \quad |\sh_{f,e}(x) - \sh_{f,a}(x)| \leq 6\nu.
\]
We use (15.4) to compare defect and maximal guaranteed vertex shift towards
$a\in M$ and $e\in \partial M$.  Choosing $\nu$ sufficiently small compared
with $\varepsilon$ (as given in (iii) of Theorem 14.5) we find -- details left
to the reader --
\[
|\mu_e({\cal F}\mid \sigma) - \mu_a({\cal F}\mid \sigma)| <
\frac{\varepsilon}{3}
\]
and
\[
|\delta_e({\cal F}\mid \sigma) - \delta_a({\cal F}\mid \sigma)| <
\frac{\varepsilon}{3}
\]
and therefore
\[
|\mu_a({\cal F}\mid \sigma) - \delta_a({\cal F}\mid \sigma)| \geq
\frac{\varepsilon}{3}
\]
for all cells $\sigma$ of $X^n$ with $h(C(\sigma)) \cap B_R(a) = \emptyset$.  By
Theorem 6.8 this establishes that $X$ is uniformly $CC^{n-1}$ over $M$.

\hfill $\Box$
~\\

We remark that Theorem 15.3 for $n=0$ recovers Theorem 12.2 (under the
assumption that $M$ is almost geodesically complete).

\newpage

\noindent {\large{\bf References}}
~\\

\fontsize{10}{10pt} \selectfont  
\bigspitem{\AA 86} H. {\AA}berg, 
Bieri-Strebel valuations (of finite rank), {\em Proc.
London Math. Soc.} (3) {\bf 52} (1986), 269-304.

\bigspitem{BeBr 97} M. Bestvina and N. Brady, Morse theory and
finiteness properties of groups, {\em Invent. Math.} {\bf 129} (1997),
445-470.

\bigspitem{B 81} R. Bieri, {\em Homological dimension of discrete groups},
QMC Mathematics Notes, 2nc edition. Queen Mary and Westfield College,
London.

\bigspitem{B 98} R. Bieri, Finiteness length and connectivity length for
groups, {\em Geometric group theory down under}, de Gruyter Verlag, Berlin
(to appear).

\bigspitem{BG$_{\mbox{I}}$} R. Bieri and R. Geoghegan, Connectivity
properties of group actions on non-positively curved spaces I: controlled
connectivity and opennenss results, preprint.

\bigspitem{BG 98} R. Bieri and R. Geoghegan, Kernels of actions on
non-positively curved spaces, {\em Geometry and cohomology in group theory}
(P.H.  Kropholler, G. Niblo and R. St{\"o}hr, ed.) London Math. Soc. 
Lecture Notes 252, Cambridge University Press, Cambridge, 1998, 24-38.

\bigspitem{BGr 84} R. Bieri and J.R.J. Groves, The geometry of the set of
characters induced by valuations, {\em J. reine und angew. Math.} {\bf 347}
(1984), 168-195.

\bigspitem{BNS 87} R. Bieri, W. Neumann and R. Strebel, A geometric
invariant of discrete groups, {\em Invent. Math.} {\bf 90} (1987), 451-477.

\bigspitem{BRe 88} R. Bieri and B. Renz, Valuations on free resolutions and
higher geometric invariants of groups, {\em Comment Math. Helvetici} {\bf
63} (1988), 464-497.

\bigspitem{BS 78} R. Bieri and R. Strebel, Almost finitely presented soluble
groups, {\em Comment. Math. Helv.} {\bf 53} (1978), 258-278.

\bigspitem{BS 80} R. Bieri and R. Strebel, Valuations and finitely presented
metabelian groups, {\em Proc. London Math. Soc.} (3) {\bf 41} (1980),
439-464.

\bigspitem{BS 81}  R. Bieri and R. Strebel, A geometric invariant for
modules over an Abelian group, {\em J. Reine und angew. Math.} {\bf 322}
(1981), 170-189.

\bigspitem{BrHa} M. Bridson and A. Haefliger, {\em Metric spaces of 
non-positive
curvature}, (monograph in preparation).

\bigspitem{Br 87$_{\mbox{II}}$} K.S. Brown, Trees, valuations and the
Bieri-Neumann-Strebel invariant, {\em Invent. Math.} {\bf 90} (1987),
479-504.

\bigspitem{BuGo 99} K.-U. Bux and C. Gonzalez, The Bestvina-Brady construction
revisited - the geometric computation of the $\Sigma$-invariants for right
angled Artin groups, {\em J. London Math. Soc.} 1999 (to appear).

\bigspitem{Ho 97} P.K. Hotchkiss, The boundary of a Busemann space, {\em
Proc. Amer. Math. Soc.} 125 (1997) 1903-1912.

\bigspitem{Ko} D. Kochloukova, {\em The $FP_m$-conjecture for a class of
metabelian groups and related topics}, Dissertation, Cambridge, 1998.

\bigspitem{Ko 96} D. Kochloukova, The FP$_m$-conjecture for a class of
metabelian groups, {\em J. of Algebra} {\bf 184} (1996), 1175-1204.

\bigspitem{Me 95} H. Meinert, The Bieri-Neumann-Strebel invariant for graph
products of groups, {\em J. Pure Appl. Algebra} {\bf 103} (1995), 205-210.

\bigspitem{Me 96} H. Meinert, The homological invariants of metabelian
groups of finite Pr{\"u}fer rank:  a proof of the $\Sigma^m$-conjecture,
{\em Proc. London Math. Soc.} (3) {\bf 72} (1996), 385-424.

\bigspitem{MMV 98} J. Meier, H. Meinert, and L. VanWyk, Higher generation 
subgroup sets and the $\Sigma$-invariants of graph groups, {\em Comment.
Math. Helv.} {\bf 73} (1998), 22-44.

\bigspitem{On} P. Ontaneda, Cocompact CAT(0) spaces are almost extendible,
(preprint).  

\bigspitem{Re 88} B. Renz, {\em Geometrische Invarianten und
Endlichkeitseigenschaften von Gruppen}, Dissertation, Frankfurt (1988).

\bigspitem{Re 89} B. Renz, Geometric invariants and HNN-extensions, {\em Group
Theory (Singapore 1987)}, de Gruyter Verlag, Berlin 1989, 465-484.
~\\

\noindent \textsc{Robert Bieri, Fachbereich Mathematik, Robert-Mayer-Strasse 
6-10, 60325 Frankfurt, Germany}

\noindent {\em e-mail}:  bieri@math.uni-frankfurt.de
~\\

\noindent \textsc{Ross Geoghegan, Department of Mathematics, State 
University of New York, Binghamton, NY  13902-6000, USA}

\noindent {\em e-mail}:  ross@math.binghamton.edu

\end{document}